\documentclass[12pt]{article}
\usepackage{latexsym}
\parskip=\medskipamount
\renewcommand{\labelenumi}{(\roman{enumi})}
\newtheorem{definition}{Definition}[section]

\newtheorem{theorem}[definition]{Theorem}
\newtheorem{proposition}[definition]{Proposition}
\newtheorem{corollary}[definition]{Corollary}
\newtheorem{remark}[definition]{Remark}
\newtheorem{example}[definition]{Example}

\font\fr=eufm10  scaled \magstep 1   
\font\ddpp=msbm10  scaled \magstep 1  

\def\QED{\hskip0.1em\hfill\null\ \null\nobreak\hfill
\kern3pt\lower1.8pt\vbox{\hrule\hbox
{\vrule\kern1pt\vbox{\kern1.7pt r\hbox{$\scriptstyle
QED$}\kern0.2pt}\kern1pt\vrule}\hrule}}
\def\R{\hbox{\ddpp R}}               
\def\C{\hbox{\ddpp C}}
\def\P{\hbox{\ddpp P}}
\def\ene{\hbox{\ddpp N}}    

\newcommand{\GL}[1]{Gl({#1},\R)}

\newcommand{\FM}{{\mathcal{F}}M}
\usepackage{amssymb}
\newcommand{\Hom}{\rm{Hom}\,}

\begin{document}

\title{\bf CLASSIFICATION OF MATERIAL $G$-STRUCTURES
\thanks{{\it (2000) MS Classification}: 53C10, 74B20, 74E05.
{\it PACS numbers}: 03.40.Dz, 02.40.Ma.
{\it Key words:} Uniformity, homogeneity, $G$-structures,
integrability, Lie groupoids, connections.
}}
\author{DAVID MAR\'IN \and MANUEL DE LE\'ON}
\date{January 25, 2004}

\maketitle

\begin{quote}
ABSTRACT. A geometrical interpretation of the $G$-structures
associated to elastic material bodies is given. In addition,
characterizations of their integrability are obtained. Since the lack
of integrability is a geometrical measure of the lack of homogeneity,
the corresponding inhomogeneity conditions are obtained.
\end{quote}

\section{Introduction}

The continuous theories of inhomogeneities were introduced by W.
Noll \cite{Noll}. In fact, Noll defined the notion of uniformity
of a hyperelastic material body using only the constitutive
law, which expresses the mechanical response of the elastic body
in terms of the gradient of the deformation. Thus, a body is
uniform if we can connect two arbitrary different points via a
material isomorphism, that is, a linear isomorphism between the
corresponding tangent spaces such that the mechanical response  at
both points is the same. The notion of material symmetry at a
point also appears in a very natural way as a linear
transformation of the tangent space at the point which does not
change the mechanical response. These notions can be translated in
a modern geometrical language in terms of Lie groupoids and Lie
groups. Indeed, the uniformity permits to construct a
$G$-structure on the body manifold whose integrability is
equivalent to the local homogeneity of the material body.

The work by Noll was extended by C.C. Wang \cite{Wang} in a
setting of principal bundles, but without an explicit mention of
the theory of $G$-structures (see also \cite{Bloom}). The first
time that the theory of $G$-structures appears explicitly linked
to uniformity occurs in a paper by Elzanowski, Epstein \&
Sniatycki \cite{EES}. In that paper, the authors have also
considered several types of $G$-structures corresponding to
different kinds of materials. However, a systematic study of the
integrability of the so-called material $G$-structures is not
available up to our knowledge. This is just the aim of the present
paper. For a material $G$-structure we mean a $G$-structure on a
material body where $G$ is a Lie subgroup of the special general
group $Sl \, (3, \R)$. We use a classification of these subgroups
usually attributed to S. Lie \cite{Lie,Nono,Wang}.

The first remarkable fact is the difficulty to obtain
integrability conditions for some of these $G$-structures in
contrast with the low dimension that we are considering. The
second point to remark is the additional difficulties arising from
the fact that we are considering subgroups of $Sl \, (3, \R)$
instead of subgroups of $Gl \, (3, \R)$. All these difficulties
are conveniently discussed along the paper.

The paper is structured as follows. In Section 2 we discuss
$G$-structures defined by tensors in a general setting (the
manifold is not necessarily three-dimensional). We give an slight
generalization of some results contained in \cite{Fu} for
nonlinear ``tensors''. Moreover, we establish some properties
concerning $G$-structures obtained by intersecting and enlarging.
These results will be very useful later. The integrability of
general $G$-structures is studied in Section 3. We propose a new
method to do this, using local $G$-connections instead of global
ones. The method leads us to integrability conditions involving
linear partial differential equations whereas the usual procedures
lead to more complicated PDE's. Section 4 is devoted to discuss
some $G$-structures defined by tensors, in particular, vector
fields, one-forms, two-forms, metrics and tensor fields of type
(1,1). We use the results previously obtained by E.T. Kobayashi
\cite{Kob2,Kob3} and J. Lehmann-Lejeune \cite{LL} for 0-deformable
tensor fields. We notice the amazing similarity between the
definition of 0-deformability in \cite{Greub} and the notion of
uniformity. In Section 5 we recall the formulation of the
continuous theories of inhomogeneities in geometrical terms. Thus,
the uniformity of the body permits to associate with it a Lie
groupoid, in such a way that, fixing a linear frame at a point (a
reference crystal) one obtains a $G$-structure, where $G$ is the
isotropy group at that point. Notice that this $G$-structure is
defined modulo conjugation, but this is sufficient for our
purposes, since the integrability is not affected by conjugation.
In Section 6, after recalling the classification of the connected
subgroups of $Sl \, (3, \R)$ modulo conjugation, we give a
geometrical interpretation of the corresponding $G$-structures,
and  we simultaneously obtain in many cases the integrability
condition. When the integrability condition is expressed in terms
of the vanishing of some tensor fields, they would be just the
inhomogeneity tensors for the corresponding material. Finally, in
Section 7, we recall a classic theorem due to Chevalley and we
give some applications. Using the natural representation, it
implies that for each algebraic subgroup $G$ of $\GL{n}$, every
${\cal N}(G)$-structure is given by the projectivization of a
tensor field which is sum of 0-deformable tensor fields, where
${\cal N }(G)$ is the normalizer of $G$ in $\GL{n}$.

\section{$G$-structures defined by tensors}

Along this paper, $\{e_{1}, \dots , e_{n}\}$ will denote the
canonical basis of $\R^{n}$, and $\{e^{1}, \dots, e^{n}\}$ its
dual basis. The space of tensors of type $(r,s)$ will be denoted
by $T_{r}^{s} \R^{n}=(\R^n)^{\otimes r}\otimes ((\R^n)^*)^{\otimes
s}$. We also notice that the action of $\GL{n}$ over ${\rm
End}(\R^{n}) = \R^{n} \otimes (\R^{n})^*$ is the functorial action
induced by the adjoint representation $A \mapsto (B \mapsto
ABA^{-1})$.

Let $M$ be an $n$-dimensional manifold and denote by $\FM$ its
linear frame bundle. $\FM$ is a principal bundle over $M$ with
projection $\pi : \FM \longrightarrow M$ and structure group
$\GL{n}$. A $G$-structure $P$ on $M$ is just a $G$-reduction of
$\FM$ (see \cite{Bernard,Chern,CDL,Fu,Kob}).

Assume that $\GL{n}$ acts on a manifold $F$ on the left. Fixing an
element $u \in F$ we denote by $G_{u}$ the isotropy group at
$u$, and by $F_{u}$ the orbit of the action through $u$. Thus,
we have
\begin{eqnarray*}
G_{u} & = & \{ a \in \GL{n} ~ | ~ au = u \} \; , \\
F_{u} & = & \{ au ~ | ~ a \in \GL{n} \} \; .
\end{eqnarray*}

\begin{definition}
An {\sf $F$-tensor} on $\FM$ is a differentiable mapping $t : \FM
\longrightarrow F$ such that $t(za) = a^{-1} t(z)$, for all $a
\in \GL{n}$ and $z \in \FM$.
\end{definition}

The following result gives the family of $G$-structures defined
by tensors. It is a slight generalization of that proved in
\cite{Fu}.

\begin{theorem}
\label{Teorema}
Giving a $G_{u}$-structure on $M$ is the same as giving an
$F$-tensor on $\FM$ which satisfy the following two conditions:
\begin{enumerate}
\item $t$ takes values in $F_{u}$;
\item $t$ is a differentiable map of $\FM$ into $F_{u}$.
\end{enumerate}
\end{theorem}

The proof is omitted, since it is a direct translation of that in
\cite{Fu}. We only remark that the relation between $t$ and the
$G_{u}$-structure is given by the formula
$P_{u} = t^{-1}(u)$.

\begin{remark}{\rm
In the case where $F_{u}$ is an embedded submanifold of $F$,
then an $F$-tensor $t$ which take values in $F_{u}$ is
automatically differentiable as a map $t : \FM \longrightarrow
F_{u}$. This is the case if $F_{u}$ is locally compact. For
instance, if $G$ is a real algebraic subgroup of $\GL{n}$ (see
\cite{Fu}).}
\end{remark}

\begin{remark}{\rm
If $u_{1}$ and $u_{2}$ are in the same orbit, say $u_{1} = a
u_{2}$, for some $a \in \GL{n}$, then $G_{u_{1}} = a G_{u_{2}}
a^{-1}$ and $P_{u_{1}} = P_{u_{2}}
a$, that is, $P_{u_{1}}$ and $P_{u_{2}}$ are conjugate.
}
\end{remark}

Since $\GL{n}$ acts on $F$ we can construct the associated
bundle $E = (\FM \times F)/\GL{n}$ over $M$ with typical fiber
$F$. Let us recall that $E$ consists of the equivalence classes
of pairs $(p,\xi) \in \FM \times F$ such that $(p, \xi) \sim
(pa, a^{-1}\xi)$.

\begin{proposition}(see \cite{Fu} for the linear case).
\label{prop}
There exists a one-to-one correspondence between $F$-tensors and
sections of $E$.
\end{proposition}

{\bf Proof:} In fact, given an $F$-tensor $t$ we define
$\sigma_{t} : M \longrightarrow E$ by $\sigma_{t}(x) = [(p,
t(p)]$, where $p$ is a linear frame at $x$.
\hfill $\Box$

\begin{remark}{\rm
The above correspondence is nothing but the extension of the
classical definition of tensor fields. Given a basis $p$ of the
tangent space $T_{x}M$ we associate
the components $(t(p))$ to the tensor $t$, which change according to
the well-known rule. It should be remarked that, if $p \in P$, where $P$
is the $G_{u}$-structure defined by $t$, then $\sigma_{t}(x) =
[(p,u)]$ (see \cite{KN}).
}
\end{remark}

\begin{corollary}
\label{cor}
Assume that $H$ is a closed subgroup of $\GL{n}$. Then there
exists a one-to-one correspondence between $H$-structures and
sections of the principal bundle $\FM/H$.
\end{corollary}

{\bf Proof:} Take $F = \GL{n}/H$ and $u = [e]$, where $e$
denotes the neutral element of $\GL{n}$. Thus, $\GL{n}$ acts on
the homogeneous space $F$ in the obvious manner, and we have
$F_{u} = F$. Since $\FM/H \cong (\FM \times \GL{n}/H)/\GL{n}$ we
deduce the result.
\hfill $\Box$

\begin{remark}
\label{rem}
{\rm
It should be noticed that Theorem \ref{Teorema}, Proposition
\ref{prop} and Corollary \ref{cor} are still true for arbitrary
principal bundles, with the obvious extension of the notion of
$F$-tensor. In this way, every $H$-reduction of the structure
group $G$ of a principal bundle $P$ to a closed subgroup $H$ may
be viewed as defined by an $G/H$-tensor on $P$.}
\end{remark}

\bigskip

Let $F$ and $F'$ be two manifolds on which $\GL{n}$ acts on the
left, and $\phi : F \longrightarrow F'$ a $\GL{n}$-invariant
differentiable mapping, i.e., $\phi(a\xi) = a \phi(\xi)$, for
all $x \in F, a \in \GL{n}$. Given a point $u \in F$ we denote
by $G_{u}$ and $G'_{u'}$ the isotropy groups of $u$ and
$u'=\phi(u)$, respectively. It si easy to check that $\phi$
induces a mapping between the associated fiber bundles, namely
$\Phi : E = (\FM \times F)/\GL{n} \longrightarrow E' = (\FM
\times F')/\GL{n}$ given by $\Phi ([p, a]) = [p, \phi(a)]$.
Thus, given a section $\sigma$ of $E$ we obtain a section $\Phi
\circ \sigma$ of $E'$. Therefore, we have obtained a way to
relate $G_{u}$ and $G'_{u'}$-structures. If $P$ is a
$G_{u}$-structure defined by an $F$-tensor $t$, we obtain a
$G'_{u'}$-structure defined by an $F'$-tensor $t'$ according to
Proposition \ref{prop}. In fact, $t$ induces a section
$\sigma_{t}$ of $E$ and $t'$ is given by the section $\Phi \circ
\sigma_{t}$. Notice that $G_{u} \subset G'_{u'}$. The
above procedure corresponds to enlarge the structure
group. Conversely, given a $G'_{u'}$-structure, we can detect
its reducibility to a $G_{u}$-structure by checking if the section
$\sigma_{t'}$ factorizes through $E$.

\begin{example}
{\rm
Consider the natural action of $\GL{n}$ on $\R^{n}$ and denote
by $S_{k}(\R^{n})$ the Stiefel manifold of $k$ frames of
$\R^{n}$. By $G_{k}(\R^{n})$ we will denote the Grassmannian of
$k$ planes in $\R^{n}$. There exists a canonical mapping
$\phi : S_{k}(\R^{n}) \longrightarrow G_{k}(\R^{n})$ which
assigns to each $k$ frame $u$ the $k$-plane $u' = \langle u
\rangle$ generated by it. Let $u$ be the $k$ frame consisting of
the $k$ first elements of the standard basis of $\R^{n}$.
A direct computation shows that
\[
G_{u} = \left\{ \left( \begin{array}{cc}
I_{k} & B \\
0&C \end{array}\right) ~ | ~ A\in\GL{k},\ C \in \GL{n-k}
\right\} \; ,
\]
where $I_{k}$ denotes the identity matrix of order $k$.
Moreover, we get
\begin{equation}
\label{grupo}
G'_{u'} = \left\{ \left( \begin{array}{cc}
A & B \\
0&C \end{array}\right) ~ | ~ A\in\GL{k},\ C \in \GL{n-k}
\right\} \; .
\end{equation}
Alternatively, we can describe $G'_{u'}$ as follows:
\[
G'_{u'} = \{ a \in \GL{n} ~ | ~ \exists \Lambda\in\GL{k},
au=\Lambda u\} \; .
\]
For $k=1$, we have that $G_1(\R^{n})$ is the projective space $\P
\R^{n}$ and $G'_{u'} = \{a \in \GL{n}~|~\exists \lambda\in \R^*,
au=\lambda u\} = \{\lambda a ~ | ~\lambda\in\R^*, a\in
GL{n}\}=\R^* G_{u}$.

It should be noticed that the action of $\GL{n}$ is transitive
and hence the orbits of $u$ and $u'$ are the whole manifolds
$S_{k}(\R^{n})$ and $G_{k}(\R^{n})$, respectively.

To end this example, take a section $\sigma_{t'}$
of $E' = (\FM \times G_{k}(\R^{n}))/\GL{n}$ associated
with an $F'_{u'}$-tensor $t'$. $\sigma_{t'}$ maps each
point $x \in M$ into a $k$-plane in $T_{x}M$, or, in other words, a
$k$-dimensional distribution on $M$. Conversely, given a
$k$-dimensional distribution on $M$, we can construct the
corresponding $F'_{u'}$-tensor.
}
\end{example}
\bigskip

Next, consider two $G$-structures $P_{1}$ and $P_{2}$ defined by
an $F_{1}$-tensor $t_{1}$ and an $F_{2}$-tensor $t_{2}$,
respectively. We assume that there is a section of ${\cal F}M$ which takes
values into $P_{1} \cap P_{2}$. Here $F_{1}$ and $F_{2}$ are manifolds
on which $\GL{n}$ acts on the left. We assume that $P_{1} =
t^{-1}(u_{1})$ and $P_{2} = t^{-1}(u_{2})$, where $u_{1} \in
F_{1}$ and $u_{2} \in F_{2}$. The corresponding structure groups
are the isotropy groups $G_{u_{1}}$ and $G_{u_{2}}$. Define an
action of $\GL{n}$ on the product manifold $F = F_{1} \times
F_{2}$ in the natural way, namely $a(\xi_{1}, \xi_{2}) = (a
\xi_{1}, a\xi_{2})$. Fixing a point $u = (u_{1}, u_{2}) \in F$,
we deduce that
\[
G_{u} = G_{u_{1}} \cap G_{u_{2}} \; , \;
F_{u} \subset F_{u_{1}} \times F_{u_{2}} \; .
\]
Define now an $F$-tensor $t$ on $\FM$ by
\[
t(p) = (t_{1}(p), t_{2}(p) \; .
\]
A direct computation shows that $t$ takes values in $F_{u}$. We
assume that $t$ is smooth as a mapping from $\FM$ into $F_{u}$
(this happens if $F_{u}$ is an embedded submanifold of $F_{1} \times
F_{2}$, for instance).
Moreover, we have $t^{-1}(u) = t_{1}^{-1}(u_{1})
\cap t_{2}^{-1}(u_{2})$, from which we deduce that $t$ defines a
$(G_{u_{1}} \cap G_{u_{2}})$-structure on $M$.

Conversely, given an $F$-tensor $t$, we can recover $t_{1}$ and
$t_{2}$ by composing $t$ with the canonical projections $F
\longrightarrow F_{1}$ and $F \longrightarrow F_{2}$. Thus, we
have proved the following.

\begin{proposition}
The intersection of two $G$-structures defined by tensors is a
new $G$-struc\-ture defined by a tensor and with structure group
the intersection of both groups.
\end{proposition}

\bigskip

Finally, a direct application of Theorem \ref{Teorema} for
arbitrary principal bundles (see Remark \ref{rem}) yields the
following construction.

Let $G$ and $G_{1}$ be closed subgroups of $\GL{n}$
such that $G \subset G_{1} \subset \GL{n}$, and assume that
$G_{1}$ acts on a manifold $F$, and $G$ is the isotropy group of
$u \in F$ under this action. Notice that we suppose that only
$G_{1}$ acts on $F$, not necessarily the whole group $\GL{n}$.

\begin{proposition}
\label{prop2}
Giving a $G$-structure on $M$ is the same as giving a
$G_{1}$-structure $P_{1}$ and an $F$-tensor $t$ on $P_{1}$ such
that
\begin{enumerate}
\item $t$ takes values in $F_{u}$.
\item $t$ is a differentiable map of $P_{1}$ into $F_{u}$.
\end{enumerate}
\end{proposition}

Of course, Proposition \ref{prop2} can be applied to the
situation of a $G_{1}$-structure and a $G_{2}$-structure defined
by two tensors, by considering $G_{1} \cap G_{2} \subset G_{1}
\subset \GL{n}$, and the action of $G_{1}$ on $F = F_{1} \times
F_{2}$. Thus, the $(G_{1} \cap G_{2})$-structure is obtained by
reducing first $\FM$ to $G_{1}$, and, then, defining an
$F$-tensor on the $G_{1}$-reduction $P_{1}$.

\begin{example}
{\rm
Consider the Grassmannian manifold $F_{1} = G_{k}(\R^{n})$, and
the natural action of $\GL{n}$ on it. Let $u_{1}$ be the $k$-plane
spanned by the first $k$ elements of the standard basis of
$\R^{k}$. Thus, $G_{u_{1}}$ is given by (\ref{grupo}).
As we know, a $G_{u_{1}}$-structure is just a $k$-dimensional
distribution $D$ on $M$. Let $F$ be the vector space of
positive definite symmetric covariant tensors of order 2 on
$\R^{k}$. $G_{u_{1}}$ acts on $\R^{k}$, but this is not the case
for $\GL{n}$! Take an inner product $u$ on $\R^{k}$, say $u \in
F$. The isotropy group $G$ of $u$ is just
\begin{equation}
\label{grupo2}
G = \left\{ \left( \begin{array}{cc}
A & B \\
0&C \end{array}\right) ~ | ~ A \in O(k),\ C \in \GL{n-k}
\right\} \; .
\end{equation}
Therefore, a $G$-structure on $M$ consists in a $k$-dimensional
distribution $D$ on $M$ endowed with an inner product on each
subspace $D_{x}$, $x \in M$. In other words, if we view a
distribution on a manifold as a vector subbundle of the tangent
bundle, a $G$-structure on $M$ consists of a vector subbundle with
a fiber metric. }
\end{example}

The above construction can be extended to include more general
structures:

\begin{itemize}

\item {\it Tangent $H$-structure on a $k$-distribution.}

Assume that $G$ is the group
\[
G = \left\{ \left( \begin{array}{cc}
A & B \\
0&C \end{array}\right) ~ | ~ A \in H \subset GL(k, \R),\ C \in \GL{n-k}
\right\} \; .
\]
Thus, giving a $G$-structure is equivalent to giving a $k$-dimensional
distribution on $M$, and a ``$H$-structure'' on each vector subspace
$D_{x}$, $x \in M$. This means that if $D$ is involutive then we have a
$H$-structure on each leaf of the induced foliation.

\item {\it Transverse $H$-structure to a $k$-distribution.}

On the other hand, assume that
\[
G = \left\{ \left( \begin{array}{cc}
A & B \\
0&C \end{array}\right) ~ | ~ A \in GL(k, \R),\ C \in H \subset \GL{n-k}
\right\} \; .
\]
Now, giving a $G$-structure is equivalent to giving a $k$-dimensional
distribution on $M$, and a ``$H$-structure'' on each quotient vector
space $T_{x}M/D_{x}$, for all $x \in M$. This means that if $D$ is
involutive then we have a foliation with transverse
$H$-structure. In such a case, we say that the $G$-structure is {\sf
projectable} if there exists a local reference $\{X_{1}, \dots, X_{n}\}$
on an open subset $U$ in $M$ such that
\[
[X_{i}, X_{a}] = 0 \;, \;  1 \leq i \leq k, 1 \leq a \leq n \; .
\]
This implies that the local quotient manifold $U'/D$ admits a
$H$-structure, where $U'$ is possible smaller than $U$.

\end{itemize}

We end this section with two examples of
tangent and transverse $H$-structures.

\begin{example}
{\rm
Let $H$ be the subgroup of $Gl(2,\R)$ given by
\[
H = \left\{\left(\begin{array}{cc}a&0\\0&c\end{array}\right)\,\big|\,
a^{\beta}=c^{\alpha}\right\},
\]
with $\alpha, \beta \in \ene$. A direct inspection shows that
$H$ is the isotropy group of the tensor
\[
\underbrace{e^1\otimes\cdots\otimes e^1}_{\beta}\otimes
\underbrace{e_2\otimes\cdots\otimes e_2}_{\alpha}.
\]

1. A tangent $H$-structure on a 2-dimensional distribution  on a
3-dimensional manifold $M$ is a $G$-structure with
\[
G=\left\{\left(\begin{array}{ccc}a&b&e\\ 0 & d & f\\ 0&0&
g\end{array}\right)\in Gl(3,\R), a^{\beta}=d^{\alpha}\right\} \; ,
\]
and it is given by a 2-dimensional distribution $D$ and a tangent
tensor field $T$ of type $(\alpha,\beta)$, i.e., a section of
$\xi^{\alpha,\beta}=\xi^{\otimes\alpha}\otimes(\xi^*)^{\otimes\beta}$
where $\xi$ is the vector sub-bundle $D\to M$ of $TM\to M$.

2. A transverse $H$-structure to a 1-dimensional distribution on a
3-dimensional manifold $M$ is a $G$-structure with
\[
G=\left\{\left(\begin{array}{ccc}e&f&g\\ 0 & a & b\\
0&0&d\end{array}\right)\in Gl(3,\R), a^{\beta}=d^{\alpha}\right\} \; ,
\]
and it is given by a 1-dimensional distribution $L$ and a transverse
tensor field $T$ of type $(\alpha,\beta)$, i.e., a section of
$\xi^{\alpha,\beta}=\xi^{\otimes\alpha}\otimes(\xi^*)^{\otimes\beta}$
where $\xi$ is the quotient vector bundle $TM/L\to M$ of $TM\to M$.
}
\end{example}

\section{Integrability}

A $G$-structure $P$ on $M$ is said to be {\sf integrable} if it is
locally equivalent to the flat standard $G$-structure $\R^n\times
G\to \R^n$, where $\dim M=n$ (see \cite{Fu}). This is equivalent
to the existence of local coordinates $(x^i)$ such that the local
section $\displaystyle{(x^i)\mapsto (x^i,\frac{\partial}{\partial
x^i})}$ is adapted.

The main problem in the theory $G$-structures is to give geometric
characterizations of their integrability. For this purpose, it is
very useful the notion of a $G$-connection.

\begin{definition}
A linear connection $\nabla$ in $M$ is said to be a $G$-connection
for a $G$-structure $P$ on $M$ if the parallel transport maps
adapted frames into adapted frames.
\end{definition}

\begin{remark}{\rm
Thus, a linear connection $\nabla$ on $M$ is a $G$-connection if
its horizontal distribution is tangent to the reduced sub-bundle,
or equivalently, $\nabla$ reduces to a connection on $P$. If the
$G$-structure $P$ is defined by a tensor $t$ then $\nabla$ is a
$G$-connection if and only if $\nabla K=0$, where $K$ is the
tensor field on $M$ defined from $t$. }\end{remark}

Since the integrability problem is a local notion, we can consider
only local $G$-connections.

\begin{proposition}
A $G$-structure $P$ is integrable if and only if around each point
in $M$ there are an open neighbourhood $U$ and a locally flat
$G$-connection $\nabla$ on $U$ (i.e., $\nabla$ is torsion less and
with zero curvature).
\end{proposition}

{\bf Proof:}

$[\Rightarrow]$ Let $(x^i)$ be local coordinates such that
$\displaystyle{(\frac{\partial}{\partial x^i})}$ is an adapted
local frame of $P$. Then the connection $\nabla$ defined by
\[
\nabla_{\frac{\partial}{\partial x^i}}\frac{\partial}{\partial
x^j}=0
\]
is locally flat.

$[\Leftarrow]$ Since the curvature of $\nabla$ is zero, the
horizontal distribution defined by $\nabla$ on ${\cal F}M$ is
involutive. Let $z(p)=(X_1(p),\dots,X_n(p))$ be an adapted frame
of $P$ at $p\in M$. The leaf trough $z(p)$ of the foliation
defined by $\nabla$ is totally contained in $P$ because of
$\nabla$ is a $G$-connection. Therefore, this leaf defines a smooth
parallel local section $(X_1,\dots,X_n)$ of $P$ over a
neighbourhood of $p$. Since $\nabla$ is torsionless and
$\nabla_{X_i}X_j=0$ we obtain that $[X_i,X_j]=0$, i.e., there are
coordinates $x^1,\dots,x^n$ such that
$\displaystyle{\{X_i=\frac{\partial}{\partial x^i}}\}$ is an
adapted local frame of $P$.

\hfill $\Box$

\begin{remark}{\rm
If $\nabla$ depends smoothly on some parameters then the local
frame $\{X_1,\dots,X_n\}$ depends also smoothly on them. }
\end{remark}

We take an adapted frame $X_1,\ldots,X_n$ on a coordinate
neighborhood $U$ and define  an auxiliary linear connection
$\widetilde{\nabla}$ on $U$ by means of
\[\widetilde{\nabla}_{X_i}X_j=0, \quad \forall i,j=1,\ldots,n.\]
In other words, $\widetilde{\nabla}$ is the linear connection
defined by the local parallelism $\{X_1,\ldots,X_n\}$. It is clear
that $\widetilde{\nabla}$ is a $G$-connection (adapted to $P$).
Moreover, any linear connection on $U$ is of the form
\[\nabla=\widetilde{\nabla}+\tau,\]
where $\tau$ is a tensor field of type $(1,2)$ on $U$.
If we put $\tau(X_i,X_j)=\tau_{ij}^kX_k$, for each $i=1,\ldots,n$
we can define the maps $\tau_i:U\to\mathfrak{gl}(n,\R)$ by putting
\begin{equation}
\label{matrix-convention} \tau_i(x)=(\tau_{ij}^k(x)).
\end{equation}

\begin{proposition}\label{G-connection}
With the above notations, $\nabla$ is a $G$-connection if and only
if for all $i=1,\ldots,n$ the maps $\tau_i$ take values in the Lie
algebra $\mbox{\fr g}$ of $G$.
\end{proposition}

{\bf Proof:}

$[\Rightarrow]$ Let ${\cal T}_i^t$ be the parallel transport
operator with respect to $\nabla$ along the integral curves of
$X_i$, then
\[
\tau_{ij}^kX_k=\nabla_{X_i}X_j=\frac{d}{ds}{\cal
T}_i^{-s}(X_j(s)).
\]
If $\nabla$ is $G$-connection then ${\cal
T}_i^{-s}(X_j(s))=a_{ij}^k(s)X_k$, where the matrix
$a_i(s)=(a_{ij}^k(s))$ belongs to $G$. Therefore
$\displaystyle{\tau_{ij}^kX_k=\left(\frac{d}{ds}a_{ij}^k(s)\right)X_k}$,
and then $\displaystyle{\tau_{ij}^kX_k =
\left(\frac{d}{ds}a_{ij}^k(s) \right) X_k \in\mbox{\fr g}}$.

\smallskip

$[\Leftarrow]$ Let $c(s)$ be a curve on $U$ and $\dot{c}(s) =
c^{\ell}(s)X_{\ell}(s)$ its tangent vector, where $X_{\ell}(s) =
X_{\ell}(c(s))$. Let $Y_{\ell}(s)$ be the parallel transport along
$c$ of the frame $\{X_1(0),\dots,X_n(0)\}$. Put
$Y_{i}(s)=a_i^j(s)X_j(s)$. We will show that $A(s)=(a_i^j(s))$
belongs to $G$. Indeed, we have
\begin{eqnarray*}
0&=&\nabla_{\dot{c}(s)}Y_i(s)=c^{\ell}(s)\nabla_{X_{\ell}(s)}a_i^j(s)X_j(s)\\
&=&c^{\ell}(s)(a_i^j(s)t_{\ell j}^k(s)X_k(s)+X_{\ell}(a_i^j(s))X_j(s))\\
&=&(c^{\ell}(s)a_i^j(s)t_{\ell
j}^k(s)+c^{\ell}(s)X_{\ell}(a_i^k(s)))X_k(s),
\end{eqnarray*}
Thus, we obtain the following identity:
\[
c^{\ell}(s)(A(s)t_{\ell}(s))+\frac{d A}{ds}=0,
\]
where $t_{\ell}(s)$ is the matrix with entries $t_{{\ell}i}^k(s)$.
Therefore we get
\[
A(s)=\exp\left(-\int_0^sc^{\ell}(u)t_{\ell}(u)du\right).
\]
We notice that the assumption $t_{\ell}(u)\in\mbox{\fr g}$ implies
$\int_0^sc^{\ell}(u)t_{\ell}(u)du\in\mbox{\fr g}$, thus, $A(s)\in
G$. Indeed, $A(s)$ belongs to connected component of the identity
of $G$. \hfill $\Box$

According to Proposition~\ref{G-connection}, if
$\nabla=\widetilde{\nabla}+\tau$ is a $G$-connection then we can
think the tensor field $\tau$ as a map
$\tau:U\to\Hom(\R^n,\mathfrak{g})$. Using this notation, the
following result will be useful in order to characterize those
$G$-connections $\nabla$ which are torsion free.
\begin{proposition}\label{Torsion}
The torsion tensor of the $G$-connection
$\nabla=\widetilde{\nabla}+\tau$ is given by
$T^{\nabla}=T^{\widetilde{\nabla}}+\partial_{\mathfrak{g}}\tau$,
where
$\partial_{\mathfrak{g}}:\Hom(\R^n,\mathfrak{g})\to\Hom(\bigwedge^2\R^n
,\R^n)$ is the map defined by
\[(\partial_{\mathfrak{g}}\tau)(x\wedge y)=\tau(x)(y)-\tau(y)(x).\]
\end{proposition}
It should be noticed that $\partial_{\mathfrak{g}}$ is just the
operator defining the Spencer cohomology of the Lie algebra
$\mathfrak{g}$ (see \cite{Fu}). In fact, after the identification
of $T_xM$ with $\R^n$ using the basis $X_1(x),\ldots,X_n(x)$,
we obtain
\[\partial_{\mathfrak{g}}\tau_{ij}^k(x)=\tau_{ij}^k(x)-\tau_{ji}^k(x).\]
Thus, the equivalence class of $T^{\widetilde{\nabla}}$ in
$\Hom(\bigwedge^2\R^n,\R^n)/\mathrm{Im}\,\partial_{\mathfrak{g}}$,
is just the first structure tensor of the $G$-structure and the
kernel of $\partial_{\mathfrak{g}}$ is the first prolongation
$\mathfrak{g}^{(1)}$ of the Lie algebra $\mathfrak{g}$.

As a consequence of Proposition~\ref{Torsion} we deduce that if
$T^{\widetilde{\nabla}}$ does not take values in
$\mathrm{Im}\,\partial_{\mathfrak{g}}$ then the $G$-structure $P$
is not integrable. Indeed, if $\nabla$ would be a free torsion
$G$-connection then
\[T^{\widetilde{\nabla}}(X_i,X_j)=-[X_i,X_j]=-\gamma_{ij}^kX_k=-(\partial_{\mathfrak{g}}\tau)_{ij}^kX_k.\]

\begin{remark}\label{Spencer-operator}
{\rm ~

\begin{itemize}
\item[(i)] If $\partial_{\mbox{\fr g}}$ is surjective then the structure
tensor is not an obstruction for the integrability.
\item[(ii)] If $\partial_{\mbox{\fr g}}$ is injective, then there is at most
one tensor $\tau$ such that $\nabla$ is a free torsion
$G$-connection. Moreover, any $G$-connection is determined by its
torsion tensor.
\item[(iii)] If $\partial_{\mbox{\fr g}}$ is bijective then there is a unique
(global) free torsion $G$-connection.
\end{itemize}
}
\end{remark}

The next issue is to investigate if, in addition to the
torsionless condition, we can choose $\tau$ in such a way that the
curvature $R^{\nabla}$ of $\nabla$ vanishes identically. In order
to do this, we try to modify $\nabla$ by adding a new tensor field
$S:U\to\Hom(\R^n,\mathfrak{g})$ such that
$\overline{\nabla}=\nabla+S$ verifies both conditions. Concerning
the first one, we have that
\[T^{\overline{\nabla}}=0\Longleftrightarrow\partial_{\mathfrak{g}}S=0\Longleftrightarrow
S\in\ker\partial_{\mathfrak{g}}=\mathfrak{g}^{(1)}.\] On the other
hand, recalling that
$\overline{\nabla}=\nabla+S=\widetilde{\nabla}+\tau+S$ and putting
$\Delta=\tau+S$, 
we obtain the following result:
\begin{proposition}\label{prop-curvature}
The curvature of $\overline{\nabla}=\widetilde{\nabla}+\Delta$
vanishes if and only if
\begin{equation}
\label{Curvature}
X_i(\Delta_j)-X_j(\Delta_i)+[\Delta_i,\Delta_j]-\gamma_{ij}^k\Delta_k=0,\quad
1\le i<j\le n,
\end{equation}
where the functions $\gamma_{ij}^k$ are defined by
$[X_i,X_j]=\gamma_{ij}^kX_k$.
\end{proposition}

{\bf Proof:} Using that
$\overline{\nabla}_{X_i}X_j=\Delta_{ij}^kX_k$ we obtain
\[\begin{array}{rclcl}
\overline{\nabla}_{X_i}\overline{\nabla}_{X_j}X_k&=&\Delta_{jk}^{\ell}\overline{\nabla}_{X_i}X_{\ell}+X_i(\Delta_{jk}^{m})X_{m}&=&\Delta_{jk}^{\ell}\Delta_{i\ell}^mX_m+X_i(\Delta_{jk}^{m})X_{m}\\
\overline{\nabla}_{X_j}\overline{\nabla}_{X_i}X_k&=&\Delta_{ik}^{\ell}\overline{\nabla}_{X_j}X_{\ell}+X_j(\Delta_{ik}^{m})X_{m}&=&\Delta_{ik}^{\ell}\Delta_{j\ell}^mX_m+X_j(\Delta_{ik}^{m})X_{m}\\
\overline{\nabla}_{[X_i,X_j]}X_k&=&\overline{\nabla}_{\gamma_{ij}^{\ell}X_{\ell}}X_k=\gamma_{ij}^{\ell}\overline{\nabla}_{X_{\ell}}X_k&=&\gamma_{ij}^{\ell}\Delta_{\ell
k}^mX_m,
\end{array}\]
so that
\[R^{\overline{\nabla}}(X_i,X_j)(X_k)=\Big(X_i(\Delta_{jk}^m)-X_j(\Delta_{ik}^m)+\Delta_{jk}^{\ell}\Delta_{i\ell}^m-\Delta_{ik}^{\ell}\Delta_{j\ell}^m-\gamma_{ij}^{\ell}\Delta_{\ell k}^m\Big)X_m.\]
Therefore, using the matrix convention introduced in
(\ref{matrix-convention}), we can write the condition
$R^{\overline{\nabla}}=0$ as the system (\ref{Curvature}). \hfill
$\Box$

\begin{remark}{\rm
It should be noticed that the above conditions involve linear
partial differential equations whereas the usual procedures lead
to more complicated PDE's.}
\end{remark}

As an application of this method we will give a characterization
of the integrability of the tangent $G$-structures.

Assume $G$ is the group
\[
G = \left\{ \left( \begin{array}{cc}
A & B \\
0&C \end{array}\right) ~ | ~ A \in H \subset GL(k, \R),\ C \in
\GL{n-k} \right\}
\]
and let $P$ be a $G$-structure. As we have seen in section 2, $P$
is the reduction of a bigger geometrical structure which consists
uniquely in a $k$-dimensional distribution $D$. If $P$ is
integrable then $D$ is involutive and it defines a foliation.
Furthermore, on each leaf $S_c$ we have an $H$-structure which is
also integrable.

Conversely, assume that the distribution $D$ is involutive and
that the induced $H$-structure on each leaf $S_c$ is integrable
(here $c$ denotes a transverse coordinate). In this case, taking a
local adapted frame $X_1,\ldots,X_n$ we can find a map
\[
\tau_c:S_c\to\rm{Hom}\,(\R^k,\mbox{\fr h})
\]
such that the connection $\nabla$ defined on $S_c$ by
\[
\nabla_{X_i}X_j = (\tau_c)_{ij}^{\ell}X_{\ell}\, ,\quad
i,j=1,\dots,k,
\]
is a locally flat $H$-connection.

\begin{proposition}\label{prop3.8}
If in addition $\tau_c$ depends smoothly on $c$, then $P$ is
integrable.
\end{proposition}

{\bf Proof:} Consider the inclusion $\rho : \mbox{\fr h}
\hookrightarrow \mbox{\fr g}$. We define a smooth map
\[
\tau : U \to \rm{Hom}(\R^n,\mbox{\fr g})
\]
by putting
\[
\tau(c,y)_{ij}^{\ell} =
\left\{\begin{array}{lcc}\rho(\tau_c(y))_{ij}^{\ell}&\rm{if}&i\le k\\
0 &\rm{if}&i>k\end{array}\right.
\]
with $j,\ell=1,\dots,n$. The connection $\nabla$ defined by $\tau$
is a $G$-connection such that its restriction to each leaf $S_c$
is locally flat. Let $\Sigma$ be a smooth transverse section to
the foliation induced by $D$ (for instance, we can parametrize
$\Sigma$ by the transverse coordinate $c$). Consider the
restriction $X_{|\Sigma}$ of the adapted frame
$X=(X_1,\ldots,X_n)$ to $\Sigma$. Since the restriction of
$\nabla$ to the leaves is flat, we can extend $X_{|\Sigma}$ to a
local adapted frame
$\overline{X}=(\overline{X}_1,\ldots,\overline{X}_n)$ defined on
$U$ by parallel transport along the leaves, so that
\[
\nabla_{\overline{X}_i}\overline{X}_j=0, \; \; {\rm if} \; 1\le
i,j\le k.
\]
Using Frobenius's Theorem we get coordinates $x^1,\dots,x^n$
such that $\displaystyle{\overline{X}_i=\frac{\partial}{\partial
x^i}}$ for $i=1,\dots,k$. \hfill $\Box$

Since the uniqueness of $\tau_c$ implies that it varies smoothly
with $c$, the following result is a consequence of the second
point of Remark~\ref{Spencer-operator}.
\begin{corollary}\label{cor3.9}
If $\partial_{\mbox{\fr h}}$ is injective then $P$ is integrable
if and only if the distribution $D$ is involutive and the
$H$-structures induced on leaves are all of them integrable.
\end{corollary}

We will treat the case in which $H$ is the following subgroup of
$\GL{2}$:
\[
H=\left\{\left(\begin{array}{cc}a&b\\0&c\end{array}\right)\,\big|\,
a^{\beta}=c^{\alpha}\right\},
\]
with $\alpha \beta \ne 0$ (see Example 2.13). This case will be
very useful in the following section.

\begin{proposition}\label{prop3.10}
With the above notations a $G$-structure is integrable if and only
if the distribution $D$ is involutive.
\end{proposition}

{\bf Proof:} We will show that every $H$-structure is integrable.
Let $\mbox{\fr h}$ be the Lie algebra of $H$:
\[
\mbox{\fr h}=\left\{\left(\begin{array}{cc}\alpha a&b\\ 0&\beta
a\end{array}\right)\,\big|\,a,b\in\R\right\}.
\]
An easy calculation show that $\partial_{\mbox{\fr h}}$ is
surjective and

\begin{equation}
\label{H1} \mbox{\fr h}^{(1)} = \{(S_1,S_2)\,|\,S_1=b_1B, \ S_2 =
\frac{b_1}{\alpha}A + b_2 B \},
\end{equation}
where
\[
A = \left(\begin{array}{cc} \alpha&0\\
0&\beta\end{array}\right),\quad
B = \left(\begin{array}{cc}0&1\\
0&0\end{array}\right)
\]
is a basis of $\mbox{\fr h}$. Let $X_1,X_2$ be an adapted local
frame of the $H$-structure, and put
$[X_1,X_2]=\gamma_1X_1+\gamma_2X_2$. Since $\partial_{\mbox{\fr
h}}$ is surjective we find $\tau=(\tau_{ij}^k)$ such that
$\partial_{\mbox{\fr h}}\tau=(\gamma_1,\gamma_2)$. For instance,
we can take
\[
\tau_1=\frac{\gamma_2}{\beta}A+\gamma_1B,\quad \tau_2=0.
\]
Taking into account (\ref{H1}) and putting
\[\left\{\begin{array}{rcl}
\Delta_1&=&\tau_1+S_1=\frac{\gamma_2}{\beta}A+(b_1+\gamma_1)B\\
\Delta_2&=&\tau_2+S_2=\frac{b_1}{\alpha}A+b_2B,
\end{array}\right.\]
the matrix differential equation
\[X_1(\Delta_2)-X_2(\Delta_1)+[\Delta_1,\Delta_2]
-\gamma_1\Delta_1-\gamma_2\Delta_2=0\] is equivalent to the system
\[
\left.
\begin{array}{rcl}
X_1 \left(\frac{b_1}{\alpha}\right) - X_2
\left(\frac{\gamma_2}{\beta} \right) -
\frac{\gamma_1\gamma_2}{\beta}-\frac{\gamma_2b_1}{\alpha}&=&0\\
 X_1(b_2)-X_2(b_1 + \gamma_1) + (\alpha-\beta)
\left(\frac{\gamma_2b_2}{\beta} -
\frac{(b_1+\gamma_1)b_1}{\alpha}\right)-\gamma_1(b_1+\gamma_1) -
\gamma_2b_2 &=& 0
\end{array}
\right\}
\]
where $b_1$ and $b_2$ are the unknown functions. The above system
of PDE's has always solution since the first equation does not
involve the unknown function $b_2$. The solution of this system
can be obtained by solving two ordinary differential equations.
Therefore these solutions depend smoothly on some parameters. We
conclude that all the $H$-structures are integrable, and from
Proposition~\ref{prop3.8} any $G$-structure with a tangent
$H$-structure is integrable if and only if the two-dimensional
distribution is involutive. \hfill $\Box$

Another interesting case is when $H$ is given by
\[
H=\left\{\left(\begin{array}{cc}a&0\\0&c\end{array}\right)\,\big|\,
a^{\beta}=c^{\alpha}\right\},
\]
with $\alpha\beta\ne 0$. One show that in this case
$\partial_{\mbox{\fr h}}$ is bijective and from
Corollary~\ref{cor3.9} we obtain the following.

\begin{corollary}\label{cor3.11}
A $G$-structure of this type is integrable if and only if the
distribution is involutive and the $H$-structures induced on the
leaves are all of them integrable. The latter occurs if and only
if the unique torsionless $H$-connection (that exists) on each
leaf has zero curvature.
\end{corollary}

Another useful remark for our purposes is the following.
\begin{remark}\label{Pbar}
{\rm  In several cases $G$ is the intersection of $\overline{G}$
and $Sl(n,\R)$, thus, a $G$-structure $P$ will be obtained as
intersection of a $\overline{G}$-structure $\overline{P}$ and a
$Sl(n,\R)$-structure given by a volume form $\Omega$, with the
compatibility condition that there exist on each point of $M$ an
adapted frame $v_1,\dots,v_n$ of $\overline{P}$ such that
$\Omega(v_1,\dots,v_n)=1$. }
\end{remark}
Concerning this situation we have the following result.

\begin{proposition}\label{prop3.12}
A $G$-structure $P$ is integrable if and only if the
$\overline{G}$-structure $\overline{P}$ is integrable and there
exist local coordinates $x^1,\dots,x^n$ adapted to $\overline{P}$
(i.e., $\displaystyle{(\frac{\partial}{\partial x^1},\dots,
\frac{\partial}{\partial x^n})\in\overline{P}}$) such that
$\displaystyle{{\cal L}_{\frac{\partial}{\partial x^i}}\Omega=0}$
for all $i=1,\ldots,n$.
\end{proposition}

Unfortunately, this proposition is too difficult to apply in the
form that it is stated because in order to show the integrability
of $P$ we need to find a privileged local coordinate system
adapted to $\overline{P}$. Therefore, we describe an alternative
approach characterizing the integrability of $P$. Instead of
expressing this condition in a privileged coordinate system
adapted to $\overline{P}$ as in Proposition~\ref{prop3.12}, we can
reformulate it in terms of an arbitrary coordinate system
$x^1,\ldots,x^n$ adapted to $\overline{P}$, i.e. such that
$\left(\frac{\partial}{\partial x^i}\right)\in \overline{P}$. We
write $\Omega=b(x^1,\dots,x^n)dx^1\wedge\dots\wedge dx^n$ and
consider the $\overline{G}$-connection $\widetilde{\nabla}$
defined by $\widetilde{\nabla}_{\frac{\partial}{\partial
x^i}}\frac{\partial}{\partial x^j}=0$. Any torsion free
$\overline{G}$-connection $\nabla$ can we written as
$\nabla=\widetilde{\nabla}+\tau$, where
$\tau:U\to\overline{\mathfrak{g}}^{(1)}\subset\Hom(\R^n,\mathfrak{g})$.
>From Proposition~\ref{prop-curvature}, the vanishing of the
curvature of $\nabla$ is equivalent to the system of PDE's:
\[\frac{\partial\tau_j}{\partial x^i}-\frac{\partial\tau_i}{\partial x^j}+[\tau_i,\tau_j]=0.\]
Finally, $\nabla$ is a $G$-connection if and only if
$\nabla\Omega=0$. Since
\[(\nabla_{\frac{\partial}{\partial x^i}}\Omega)(\frac{\partial}{\partial x^1},\ldots,\frac{\partial}{\partial x^n})=\frac{\partial b}{\partial x^i}-\sum_{j=1}^n\Omega(\frac{\partial}{\partial x^1},\ldots,\nabla_{\frac{\partial}{\partial x^i}}\frac{\partial}{\partial x^j},\ldots,\frac{\partial}{\partial x^n})=\frac{\partial b}{\partial x^i}-b\,{\rm tr}\,(\tau_i),\]
we can characterize the equation $\nabla\Omega=0$ as
$\frac{\partial \log(b)}{\partial x^i}={\rm tr}\,(\tau_i)$ for all
$i=1,\ldots,n$.

\begin{remark}\label{bk}
{\rm  In several cases, straightforward computations show that if
$\tau$ takes values in $\overline{\mbox{\fr g}}^{(1)}$ then ${\rm
tr}(\tau_k)=0$ for some $k$. Therefore, assuming that
$\overline{P}$ is integrable, a \emph{necessary} condition for the
integrability of $P$ is that $b(x_1,\ldots,x_n)$ does not depend
on $x^k$, or equivalently,
\[
\displaystyle{{\cal L}_{\frac{\partial}{\partial x^k}}\Omega=0}.
\]
}
\end{remark}

\section{Some examples of $G$-structures defined by tensors}

The main purpose of this paper is a systematic study of those
$G$-structures associated to uniform elastic bodies. Before to do
that, we will discuss $G$-structures defined by tensors of type
$(r,s)$, with $r+s\le 3$. The results will be useful in the next
sections.

\begin{itemize}

\item If $F=T_1^0\R^3$ then, since $\GL{3}$ acts transitively on $F$,
we can take $u=e_1$. Therefore the isotropy group $G$ of $u$ consists
of matrices of the form
\begin{equation}\label{m1}
\left(\begin{array}{ccc}1&a&b\\ 0&c&d\\ 0&e&f\end{array}\right)
\end{equation}
The associated fiber bundle is $TM$. Thus, a $G$-structure is given by
a vector field $X$ without zeros. Hence every $G$-structure
is integrable, since we can always choose local coordinates
$x^{1}, x^{2}, x^{3}$ such that
$\displaystyle{X=\frac{\partial}{\partial x^1}}$.

\item If $F=T_0^1\R^3$ then we can take $u=e^1$, and in this case $G$
is the group of matrices obtained by transposing~(\ref{m1}). The
associated fiber bundle is the cotangent bundle $T^*M$. Therefore
a $G$-structure is given by a one-form $\omega$ without zeros. Its
integrability is equivalent to the existence of local coordinates
such that $\omega=dx^1$, i.e., $\omega$ is locally exact, or
equivalently, $\omega$ is closed.

\item If $F=T_1^1 \R^{3}$, then the action of $\GL{3}$ on $F$ is by
conjugation, so that the corresponding orbits are not trivial:
$u_1,u_2\in F$ are in the same orbit if and only if they have the
same canonical form over $\R$. We will study the different
possibilities in dimension 3. The minimum polynomial is one of the
following types:
\renewcommand{\labelenumi}{(\alph{enumi})}
\begin{enumerate}
\item $((x-\alpha)^2+\beta^2)(x-\lambda)$, $\beta\ne 0$ with canonical
form over $\R$
\[
\left(\begin{array}{ccc}\alpha &\beta& 0\\ -\beta &\alpha & 0\\ 0&
0&\lambda\end{array}\right)
\]

\item $(x-\lambda)(x-\mu)(x-\nu)$, where $\lambda,\mu,\nu$ are three
different eigenvalues; in this case the Jordan form is diagonal.

\item $(x-\lambda)^2(x-\mu)$, $\lambda\ne\nu$ with Jordan form
\[
\left(\begin{array}{ccc} \lambda&1&0\\ 0&\lambda &0\\ 0&0& \mu
\end{array}\right)
\]

\item $(x-\lambda)(x-\mu)$, $\lambda\ne\nu$ with diagonal Jordan form
\[
\left(\begin{array}{ccc} \lambda&0&0\\ 0&\lambda &0\\ 0&0& \mu
\end{array}\right)
\]

\item $(x-\lambda)^3$, with Jordan form
\[
\left(\begin{array}{ccc} \lambda&1&0\\ 0&\lambda &1\\ 0&0& \lambda
\end{array}\right)
\]

\item $(x-\lambda)^2$, with Jordan form
\[
\left(\begin{array}{ccc} \lambda&1&0\\ 0&\lambda &0\\ 0&0& \lambda
\end{array}\right)
\]

\item $(x-\lambda)$, in this case we have a homothetic transformation.
\end{enumerate}\renewcommand{\labelenumi}{(\roman{enumi})}

The isotropy groups of theses matrices consist of the
matrices of the following form:

$$\begin{array}{ccc}

{\rm (a)} \left(\begin{array}{ccc} a&0&0\\ 0&b&c\\ 0&-c&b \end{array}\right)
& \hspace{10mm} &
{\rm (b)} \left(\begin{array}{ccc} a&0&0\\ 0&b&0\\ 0&0&c \end{array}\right)
\\
{\rm (c)} \left(\begin{array}{ccc} a&b&0\\ 0&a&0\\ 0&0&c \end{array}\right)
& \hspace{10mm} &
{\rm (d)} \left(\begin{array}{ccc} a&0&0\\ 0&b&c\\ 0&d&e \end{array}\right)
\\
{\rm (e)} \left(\begin{array}{ccc} a&b&c\\ 0&a&b\\ 0&0&a \end{array}\right)
& \hspace{10mm} &
{\rm (f)} \left(\begin{array}{ccc} a&0&b\\ c&d&e\\ 0&0&d \end{array}\right)
\end{array}$$

The case (g) is trivial since the homothetic transformation
commutes with every element in $\GL{3}$.

\begin{remark}
{\rm
We notice that the isotropy group is determined by
the relations between $\lambda,\mu$ and $\nu$ and the fact that
$\beta\ne 0$, but not by the particular values of
them. For instance, we can take $\lambda=0,\mu=1,\nu=-1,\alpha=0$ and
$\beta=1$. Then the non trivial orbits of $F=T_1^1M$ are
given by $u\in F$ fulfilling one and only one of the following
equations:
\[
u^3+u=0,\ u^3-u=0,\ u^3-u^2=0,\ u^2-u=0,\ u^3=0,\ u^2=0.
\]
}
\end{remark}

The associated fiber bundle is $T_1^1 M=TM\otimes T^*M$.
Therefore, a $G$-structure with $G$ a subgroup of the above list is
given by a tensor field $h$ of type (1,1) such that for each $x\in M$
$h_x:T_xM\to T_xM$ has constant canonical form $u$.
These $G$-structures were extensively treated in the literature.
The tensor field $h$ is called a 0-deformable vector
one-form (or $(1,1)$-type tensor field) \cite{LL}. The general notion
of 0-deformability can be found in \cite{Greub}:

\begin{definition}
A section $\sigma$ of a vector bundle $\xi=(E,\pi,M,F)$ is said to be
0-deformable if for each $x,y\in M$ there exist a linear isomorphism
$\alpha_{x,y}:F_x\to F_y$ such that
$\alpha_{x,y}(\sigma(x))=\sigma(y)$.
\end{definition}

\begin{remark}
{\rm
\begin{enumerate}

\item We notice the similarity of this definition with the notion of
material uniformity of a body $B$ (see Section 5).

\item A finite set of 0-deformable cross-sections
$\Sigma_{\xi}=(\sigma_1,\ldots,\sigma_m)$ of
$\xi^{p_i,q_i}=\xi^{\otimes p_i}\otimes(\xi^*)^{\otimes q_i}$ define a
$\Sigma$-bundle (see \cite{Greub}) where $\xi=(E,\pi,B,F)$ is a vector
bundle. In fact, the definition of $\Sigma$-bundle involves:

\begin{enumerate}

\item a smooth vector bundle $\xi=(E,\pi,B,F)$.

\item a finite ordered set of cross-sections
$\Sigma_{\xi}=(\sigma_1,\ldots,\sigma_m)$, $\sigma_i \in{\rm
Sec}\,\xi^{p_i,q_i}$,

\end{enumerate}

subject to the following condition:

there are a finite ordered system $\Sigma_F=(u_1,\ldots,u_m)$ of tensors
$u_i\in F^{p_i,q_i}$ and a coordinate representation
$\{U_{\alpha},\psi_{\alpha}\}$ of $\xi$ such that
\[
\psi_{\alpha}^{p_i,q_i}(x,u_i)=\sigma_i(x)\in\xi_x^{p_i,q_i},\quad
x\in U_{\alpha},\quad i=1,\ldots,m.
\]
Then we will say that the section $\sigma_i$ is 0-deformable to
$u_i$. A $\Sigma$-bundle gives a reduction of the structure group
from $GL(F)$ to an algebraic subgroup. If $\xi$ is the tangent
bundle we obtain in this way a $G$-structure with $G$ an algebraic
subgroup of $\GL{n}$.
\end{enumerate}
}
\end{remark}

Coming back to the case of a 0-deformable vector one-form $h$ we are
interested in the characterization of the integrability of the
$G$-structure defined by $h$. The following theorem was
proved by E.T. Kobayashi \cite{Kob3} (see also \cite{LL}).

\begin{theorem}
Let $h$ be a 0-deformable vector one-form on a manifold $M$, with
characteristic polynomial
\[
\prod_ip_i(x)^{d_i},
\]
where $p_i(x)$ are irreducible and coprime polynomials in $x$ over
$\R$. The minimum polynomial is
\[
\prod_ip_i(x)^{v_i}.
\]
We assume that for each $i$ we have $v_i=1$ or $v_i=d_i$.
Then the $G$-structure defined by $h$ is integrable if the Nijenhuis
tensor of $h$, $H_h$ is zero, where $N_h$ is a vector 2-form defined by
\[
\frac{1}{2}N_h(X,Y)=[hX,hY]-h[hX,Y]-h[X,hY]+h^2[X,Y] .
\]
\end{theorem}

\begin{remark}{\rm
\begin{enumerate}
\item In dimension three the above theorem gives a suficient condition
for the integrability in all the cases except when $u^2=0$.
\item $N_h=0$ is always a necessary condition for the integrability.
\end{enumerate}
}
\end{remark}

\item If $F=T_0^2\R^3=S^2\R^{3}\oplus\bigwedge^2\R^3$, the action of an
element $A\in \GL{3}$ on $f\in F$ is given by $A^t \, f\,A$. There
are two fundamental cases that we will discuss separately:

\begin{itemize}
\item[-] $u \in S^2\R^3$, then applying Sylvester's Theorem we deduce
that $u$ is in the orbit of $\epsilon_1 e^1\otimes e^1+\epsilon_2
e^2\otimes e^2+\epsilon_3 e^3\otimes e^3$ for some
$\epsilon_1,\epsilon_2,\epsilon_3\in\{-1,0,1\}$, i.e., the range
and the signature determine the orbit of $u$. The associated fiber
bundle is $T_0^2M$, and a $G$-structure with $G$ the isotropy
group of $u$ is given by a symmetric (0,2)-tensor field which is
0-deformable to $u$. If $u=e^1\otimes e^1+e^2\otimes e^2\pm
e^3\otimes e^3$ we obtain, respectively, a Riemannian or
Lorentzian metric $g$ on $M$. The integrability of this structure
is characterized by the vanishing of the scalar curvature of $g$,
but, in dimension three, this is equivalent to the vanishing of
the Ricci tensor of $g$.

\item[-] $u\in\bigwedge^2\R^3\setminus\{0\}$. By Darboux's Theorem $u$
is in the orbit of $e^1\wedge e^2$. Therefore a $G$-structure with
\[
G=\left\{\left(\begin{array}{cc}A&0\\
v&c\end{array}\right)\,\big|\,A\in Sl(2,\R),\,c\in\R^*\right\}
\]
is given by a two-form $\eta$ without zeros. The integrability is
also characterized by Darboux's Theorem:
The $G$-structure is integrable if and only if $d\eta=0$.
\end{itemize}

\item The case $F=T_2^0\R^3$ is formally analogous to the above
case, but the integrability condition is different.

\item Finally, we will consider the case $F=\bigwedge^3\R^3\subset
T_0^3\R^3$. We can take $u=e^1\wedge e^2\wedge e^3$,
$G=G_u=Sl(3,\R)$ and $Q=\Omega^3(M)$. Therefore a $G$-structure
over $M$ is given by a three-form without zeros, i.e., a volume
form $\Omega$ on $M$. As we know, every $Sl(3,\R)$-structure is
integrable.

\end{itemize}

\section{Uniformity and homogeneity of simple materials}

A {\sf body} $B$ is a 3-dimensional differentiable manifold
which can be covered with just one chart. An embedding $\Phi : B
\longrightarrow \R ^{3}$ is called a {\sf configuration} of $B$.
The body is identified with any one of its configurations, say
$\Phi_{0} : B \longrightarrow \R^{3}$, called a {\sf reference
configuration}. Given any arbitrary configuration $\Phi : B
\longrightarrow \R^{3}$, the change of configurations $\kappa =
\Phi \circ \Phi_{0}^{-1}$ is called a {\sf deformation}.
We fix a reference configuration $\Phi_{0}$ and, from now on,
$B$ and its image $\Phi_{0}(B)$ will be identified.
The mechanical behaviour of a {\sf hyperelastic} material body
is characterized by one function $W$ which depends, at each
point of $B$, only on the value of the derivative of the
deformation evaluated at that point (see
\cite{marsden1,Maugin,Noll,TrToup,TrNoll,WangTr,Wang}).
$W$ measures the strain energy per unit volume of reference
configuration. In a more general material bodies, $W$ can depend also
of higher order gradients or even more complicated microstructures (see
\cite{LE1,LE2,LE230,LE231,LE3}.

$B$ is said to be {\sf materially uniform} if for two arbitrary
points $X, Y \in B$ there exists a local diffeomorphism $\phi$
from a neighbourhood of $X$ onto a neighbourhood of $Y$ such that
$\phi(X) = Y$ and
\begin{equation}
W(j^{1}_{Y, \kappa(Y)} \kappa) = W(j^{1}_{Y,
\kappa(Y)} \kappa \cdot j^{1}_{X,Y} \phi) \; \; ,
\label{1}
\end{equation}
for all $j^{1}_{Y, \kappa(Y)} \kappa$.
The 1-jet $j^{1}_{X,Y} \phi$ will be called a {\sf material 1-jet}.
It should be noticed that $j^{1}_{X,Y} \phi$ may be identified
with the linear isomorphism $d\phi(X) : T_{X} B
\longrightarrow T_{Y}B$. $d\phi(X)$ is usually called a
{\sf material isomorphism} (\cite{TrNoll}).

Denote by $\Omega(B)$ the collection of all material
1-jets. Thus, $\Omega(B) \subset \Pi^{1}(B,B)$, where $\Pi^{1}(B, B)$
is the Lie groupoid of all invertible 1-jets on the manifold $B$ (see
\cite{Mack} for a general reference on Lie groupoids).
We have canonical mappings
$\alpha : \Pi^{1}(B, B) \longrightarrow B$
and $\beta : \Pi^{1}(B, B) \longrightarrow B$ defined by
\[
\alpha(j^{1}_{X,Y} \phi) = X \; , \;
\beta(j^{1}_{X,Y} \phi) = Y \; ,
\]
respectively. Their restrictions to $\Omega(B)$ will be
denoted by the same symbols. A direct computation from (\ref{1})
shows the following result.

\begin{proposition}
If $B$ is uniform, then $\Omega(B)$ is a groupoid
with source and target mappings $\alpha$ and $\beta$,
respectively. In fact, $\Omega(B)$ is a subgroupoid of
the Lie groupoid $\Pi^{1}(B, B)$ of all invertible
1-jets on the manifold $B$.
\end{proposition}

\begin{definition}
A {\sf material symmetry} at a point $X$ is a 1-jet $j^{1}_{X,X}\phi$
of a local diffeomorphism at $X$ such that
\begin{equation}
W(j^{1}_{X, \kappa(X)} \kappa) = W(j^{1}_{X,
\kappa(X)} \kappa \cdot j^{1}_{X,X} \phi) \; \; ,
\label{2}
\end{equation}
for all $j^{1}_{X, \kappa(X)} \kappa$.
\end{definition}

>From (\ref{2}) we deduce that the collection $G(X)$ of all
material symmetries at $X$ has a structure of group which is
called the {\sf symmetry group} at $X$. If $j^{1}_{X,Y} \phi$ is a
material 1-jet joining two points $X$ and $Y$, we deduce that the
symmetry groups at $X$ and $Y$ are conjugate:
\[
G(Y) = j^{1}_{X,Y}\phi \circ G(X) \circ (j^{1}_{X,Y} \phi)^{-1}
\; .
\]

\begin{definition}
We say that $B$ enjoys {\sf smooth uniformity} if
$\Omega(B)$ is a Lie groupoid, which will be called the
{\sf material Lie groupoid}.
\end{definition}

In such a case, there exist local sections of the projection
$(\alpha, \beta) : \Omega(B) \longrightarrow B
\times B$ given by $(\alpha, \beta)(j^{1}_{X,Y}\phi) =
(X, Y)$. Such a local section, say $P : B \times
B \longrightarrow \Omega(B)$ assigns to each pair $(X,
Y)$ of material points a material 1-jet connecting them.
Such a section $P$ is called a {\sf local material uniformity}. If
there exists a global section of $(\alpha, \beta)$, then $B$ enjoys
smooth global uniformity, or, equivalently, the Lie groupoid
$\Omega(B)$ is smoothly transitive.

By applying well-known results on Lie groupoids and frame
bundles, we get the following (see \cite{Fu}).

\begin{proposition}
\label{p2}
If $B$ enjoys smooth uniformity, then:

\begin{enumerate}
\item $G(X_{0})$ is a Lie group.
\item $\alpha^{-1}(X_{0})$ is a principal $G(X_{0})$-bundle
over $B$ whose canonical projection is the restriction of
$\beta$.
\end{enumerate}
\end{proposition}

{\bf Proof:} (i) Since $(\alpha, \beta)$ is a submersion, we deduce that
$G(X_{0})=(\alpha, \beta)^{-1}(X_{0})$ is a closed submanifold of
$\Omega(B)$. Hence, it is a Lie group.

(ii) First of all, since $\alpha$ is a submersion, we deduce that
$\alpha^{-1}(X_{0})$ is a closed submanifold of $\Omega(B)$. Moreover,
since $(\alpha, \beta)$ is a submersion, there
exist an open covering $\{U_{a}\}$ of $B$ and local
sections $\sigma_{a,b} : U_{a} \times U_{b} \longrightarrow
\Omega(B)$ of $(\alpha, \beta)$.
If $X_{0} \in U_{a_{0}}$, we define
\[
\sigma_{a}(X) = \sigma_{a_{0},a}(X_{0},X) \; , \; {\rm for} \; X
\in U_{a} \; .
\]
In other words, $\sigma_{a}$ assigns (in a differentiable way)
to each material point $X$ a material 1-jet connecting $X_{0}$
and $X$. Thus, we have obtained a family of local sections
$\{\sigma_{a}\}$ of $\beta: \alpha^{-1}(X_{0}) \longrightarrow
B$ which define a principal $G(X_{0})$-bundle.
\hfill $\Box$

Next, we fix a point $X_{0}$ at $B$. The tangent bundle
$T_{X_{0}}B$ is a linear approximation of an
infinitesimal piece of material around $X_{0}$. But
$T_{X_{0}}B$ is completely characterized by a basis. This
fact leads us to the following definition.

\begin{definition}
A linear frame $Z_{0}$ at a material point $X_{0}$ will be
called a {\sf reference crystal} (at $X_{0}$).
\end{definition}

A reference crystal $Z_{0}$ is just a 1-jet $j^{1}_{0,X_{0}}
\psi$ of a local diffeomorphism from $0 \in \R^{3}$ into
$X_{0}$. Thus, we can transport $Z_{0}$ to any point of $B$ by
composing it with smooth material uniformities.
The next result is also standard in the literature (see
\cite{Fu}).

\begin{theorem}
\begin{enumerate}
\item $G = Z_{0}^{-1} \circ G(X_{0}) \circ Z_{0}$ is a Lie
subgroup of $\GL{3}$.
\item Denote by $\omega(B)$ the set of all linear frames at
all the point of $B$ obtained by translating $Z_{0}$.
Then $\omega(B)$ is a $G$-structure on $B$.
\end{enumerate}
\end{theorem}

{\bf Proof:} (i) Put $GL(B, X_{0}) = \{
j^{1}_{X_{0},X_{0}} \phi\} \subset \Pi^{1}(B, B)\}$.
$GL(B, X_{0})$ is a Lie group which is isomorphic to
$\GL{3}$. Since the mapping
\[
GL(B, X_{0})  \longrightarrow \GL{3} \; , \;
j^{1}_{X_{0},X_{0}} \phi  \leadsto Z_{0}^{-1} \circ
j^{1}_{X_{0},X_{0}} \phi \circ Z_{0} \; ,
\]
is smooth, it follows that $G = Z_{0}^{-1} \circ G(X_{0}) \circ
Z_{0}$ is a Lie subgroup of $\GL{3}$.

(ii) Take the family of local sections $\{\sigma_{a}\}$
obtained in Proposition \ref{p2}. A family of local
sections $\{\tau_{a}\}$ of $\omega(B)$ for the open covering
$\{U_{a}\}$ is obtained as follows:
\[
\tau_{a}(X) = \sigma_{a}(X) \circ Z_{0} \; , \; {\rm for} \; X
\in U_{a} \; .
\]
An straightforward computation shows that $\omega(B)$ is
in fact a $G$-reduction of ${\cal FB}$.
\hfill $\Box$

This $G$-structure will be called {\sf material}.

\begin{remark}\label{rempi1}
{\rm
(1) If we perform a change of reference configuration, the
$G$-structure remains the same, provided that the point $X_{0}$
and the reference crystal $Z_{0}$ are dragged by the change of
configuration.

(2) If we choose a different point $X'_{0}$, we obtain the same
$G$-structure provided that the reference crystal is the one
obtained using a material uniformity from $X_{0}$ to
$X'_{0}$.

(3) If we change the reference crystal $Z_{0}$ to $Z'_{0} =
Z_{0} A$, where $A \in \GL{3}$, we obtain a conjugate
$G$-structure $\omega(B)A$, with conjugate structure
group $A^{-1}GA$.
}
\end{remark}

\begin{definition}
A body $B$ is said to be homogeneous if there exists a global
deformation $\kappa$ such that $Q$ defined by
\[
Q(X)=j^1_{0,X}(\kappa^{-1}\circ\tau_{\kappa(X)}),\;\forall X\in B,
\]
is a uniform reference, where $\tau_{\kappa(X)}:\R^3\to\R^3$ denotes
the translation by $\kappa(X)$.

$B$ is said to be locally homogeneous if around each point $X$
of $B$ there exists an open neighbourhood $U$ which is homogeneous.
\end{definition}

The following result gives a geometric characterization of the local homogeneity.

\begin{theorem}
$B$ is locally homogeneous if and only if the associated
material $G$-structure is integrable.
\end{theorem}

\begin{remark}{\rm
According to Remark~\ref{rempi1} the above definition does not depend
on the chosen crystal reference.
}
\end{remark}

\section{Classification of material $G$-structures}

Our purpose is to study systematically the possible material
$G$-structures associated to elastic bodies.

For physical reasons (see \cite{TrNoll,Wang}) we are only
interested in $G$-structures with $G$ a Lie subgroup of the
special linear group $Sl(3,\R)$. The first step is to classify the
subgroups of $Sl(3,\R)$. A classification modulo conjugation is
usually attributed to S. Lie \cite{Lie,Nono,TrNoll,Wang}. We
reproduce here the list as it is presented in Wang \cite{Wang}.
This list gives the classification of the Lie subalgebras of the
Lie algebra ${\fr sl}(3)$ of $Sl(3,\R)$ and their corresponding
connected Lie subgroups, see Appendix A.

There are three types of solids, namely:
\begin{itemize}

\item isotropic solids belong to type 16,
\item transversely isotropic solids belong to type 8 with parameter
$\alpha=0$, and

\item crystalline solids belong to type 5 with $\alpha=\beta=\gamma=0$.
\end{itemize}
All other types are fluid crystals. For instance,
\begin{itemize}
\item isotropic fluids belong to type 9, and
\item fluid crystal of first kind (respectively, second kind) belong to
type 11 (respectively, 10).
\end{itemize}

The first five families consist of three fundamental types, denoted by
A,\,B and C, respectively:
\begin{itemize}
\item type A is characterized by considering $\alpha, \beta$ and
$\gamma$ as variables and they are algebraic subgroups.
\item type B is characterized by considering $\alpha, \beta$ and
$\gamma$ as fixed parameters, with $(\alpha,\beta,\gamma)\ne(0,0,0)$
and $\alpha+\beta+\gamma=0$. They define an element
$[\alpha,\beta]$ of $\P\R^1$. The corresponding subgroup
of the list is not algebraic, but it is contained in a bigger ``natural''
algebraic subgroup if the parameters are integer.
\item type C is obtained by taking $\alpha=\beta=\gamma=0$ and all of
them are algebraic.
\end{itemize}

The types 6-8 contain two different cases:
\begin{itemize}
\item type A is characterized by taking $\alpha$ and $\beta$ as
variables. They are algebraic subgroups.
\item type B is characterized by taking $\alpha$ and $\beta$ as fixed
parameters. They are algebraic if and only if $\alpha=0$.
\end{itemize}

The other families consist of a unique type, and they are
algebraic except families 17, 18, 21, 22 and 25.

In what follows, we will discuss the $G$-structures with $G$ an
algebraic subgroup of $Sl(3,\R)$. We remark that the Lie subgroups
included in the list are connected, however we will consider the
corresponding natural algebraic subgroups.

Given a such $G$-structure, and according to Remark~\ref{Pbar},
sometimes we will consider an enlarged structure $\overline{P}$
with structural group $\overline{G}\subset\GL{3}$ such that
$\overline{G}\cap Sl(3,\R)=G$. The relation between the
integrability of $P$ and $\overline{P}$ is given by
Proposition~\ref{prop3.12}. Thus, in some cases we have only
characterized the integrability of $\overline{P}$.

\bigskip

{\bf The group 1A}

\medskip

The group $G_{1A}$ is just the isotropy group of the linear subspaces
$\langle e_3\rangle\subset\langle e_2,e_3\rangle$ and the tensor
$w=e^1\wedge e^2\wedge e^3$ on $\R^3$. Then such a $G_{1A}$-structure $P$ is
given by a one-dimensional distribution $L$, a two-dimensional
distribution $D$ with $L\subset D$, and a volume form $\Omega$.

\begin{proposition}\label{prop1A}
$P$ is integrable if and only if $D$ is involutive.
\end{proposition}

{\bf Proof:}
If $D$ is involutive and $X,Y,Z$ is an adapted local basis, i.e.,
$L=\langle X\rangle$, $D=\langle X,Y\rangle$, then we have $[X,Y]=\alpha
X+\beta Y$. Therefore there exist local functions $f$ and $g$ such that
$[fX,gY]=0$. Thus, there are local coordinates $y^1,y^2,y^3$ such that
$\displaystyle{fX=\frac{\partial}{\partial y^1}}$ and
$\displaystyle{gY=\frac{\partial}{\partial y^2}}$. If
$\Omega=b(y^1,y^2,y^3)dy^1\wedge dy^2\wedge dy^3$, we define
new coordinates $x^1=\int b(y^1,y^2,y^3)dy^1$, $x^2=y^2$, $x^3=y^3$,
and hence $(x^1,x^2,x^3)$ are adapted coordinates. The converse is
trivial.
\hfill $\Box$

\bigskip

{\bf The group 1B}

\medskip

As we have said before we only discuss the algebraic case. Thus,
$\alpha, \beta$ and $\gamma$ are integer parameters. We can only
consider the following six cases:
\begin{itemize}
\item $\alpha,\gamma>0$, then the tensor
\[
t=\underbrace{e^1\otimes\cdots\otimes e^1}_{\gamma}\otimes
\underbrace{e_3\otimes\cdots\otimes e_3}_{\alpha}
\]
is invariant with respect to the natural action of $GL(3,\R)$. In
fact, the isotropy group of $t$ and $w$ is just the group $G_{1B}$.
Therefore, in this case a $G_{1B}$-structure $P$ is given by a tensor
field $T$ of type $(\gamma,\alpha)$ which is 0-deformable to $t$, and a
volume form $\Omega$. In addition we have a two-dimensional
distribution $D$ on $B$ such that $T$ is tangent to $D$.

>From Propositions~\ref{prop3.10} and \ref{prop3.12}, we deduce the
following.

\begin{proposition}\label{prop1B}
The $\overline{G}_{1B}$-structure $\overline{P}$ is integrable if
and only if the distribution $D$ is involutive. In this case, let
us denote by $(x^1,x^2,x^3)$ local coordinates adapted to
$\overline{P}$. If $P$ is integrable then necessarily
$\displaystyle{{\cal L}_{\frac{\partial}{\partial
x^2}}\Omega={\cal L}_{\frac{\partial}{\partial x^3}}\Omega=0}$.
\end{proposition}

{\bf Proof:} It  only remains to prove the last assertion, which
follows directly from Remark~\ref{bk} and the following
computation of
$\overline{\mathfrak{g}}^{(1)}=\ker\partial_{\overline{\mathfrak{g}}}$,
where
$\partial_{\overline{\mathfrak{g}}}:\Hom(\R^3,\overline{\mathfrak{g}})\to\Hom(\bigwedge^2\R^3,\R^3)$
is given by
$\partial_{\overline{\mathfrak{g}}}((\tau_i))=(\tau_{ij}^k-\tau_{ji}^k)$
and $\tau_i=\left(\begin{array}{ccc} \alpha a_i & 0 & 0\\ b_i &
\beta a_i & 0\\ c_i & d_i & \gamma a_i
\end{array}\right)$.
>From the following table,
\begin{center}
\begin{tabular}{c|c|c|c|}
~ & $\tau_{12}^k-\tau_{21}^k$ & $\tau_{13}^k-\tau_{31}^k$ &
$\tau_{23}^k-\tau_{32}^k$\\
\hline
 $k=1$ & $0-\alpha a_2$ & $0-\alpha a_3$ & $0-0$\\
\hline
 $k=2$ & $\beta a_1-b_2$ & $0-b_3$ & $0-\beta a_3$\\
\hline
 $k=3$ & $d_1-c_2$ & $\gamma a_1-c_3$ & $\gamma a_2-d_3$\\
\hline
\end{tabular}
\end{center}
we deduce that if
$\tau=(\tau_i)_{i=1}^3\in\overline{\mathfrak{g}}^{(1)}$, then
$a_2=a_3=0$ and therefore
$\mathrm{tr}\,\tau_2=\mathrm{tr}\,\tau_3=0$. \hfill $\Box$

\item $\beta,\gamma>0$, then
\[
t=\underbrace{e^2\otimes\cdots\otimes e^2}_{\gamma}\otimes
\underbrace{e_3\otimes\cdots\otimes e_3}_{\beta}
\]
is a tensor of type $(\gamma,\beta)$ defined on the subspace
$\langle e_2,e_3\rangle$. Our group is just the isotropy group of
$t$ and $w$. A $G$-structure is now given by a two-dimensional
distribution $D$ and a tangent tensor field $D$ of type
$(\gamma,\beta)$ which is 0-deformable to $t$. The integrability
condition is the same as in the precedent case. However, it should
be noticed that now $T$ is not a global tensor field on $B$.

\item $\alpha, \beta>0$, then the subspace $\langle e_3\rangle$ and the
transverse tensor
\[
t=\underbrace{u^1\otimes\cdots\otimes
u^1}_{\beta}\otimes\underbrace{u_2\otimes\cdots\otimes u_2}_{\alpha}
\]
where $u_1,u_2$ is a basis of the quotient vector space
$\R^3/\langle e_3\rangle$, together with the three-form $w$
determine the group $G_{1B}$. Therefore a $G_{1B}$-structure $P$
is, in this case, given by a one-dimensional distribution $L$, a
transverse tensor field $T$ to $L$ which is 0-deformable to $t$,
and a volume form $\Omega$.

\begin{proposition}
$\overline{P}$ is integrable if and only if it is projectable,
i.e., the coefficients of $T$ in a system of coordinates
$x^1,x^2,x^3$ adapted to the foliation induced by $L$ does not
depend on the tangent coordinate $x^3$. Moreover, assume that
$(x_1,x_2,x_3)$ are local coordinates adapted to $\overline{P}$
and $P$ is integrable. Then $\displaystyle{{\cal
L}_{\frac{\partial}{\partial x^2}}\Omega={\cal
L}_{\frac{\partial}{\partial x^3}}\Omega=0}$.
\end{proposition}

\begin{remark}{\rm
An alternative description of $P$ is a tangent $H$-structure on a
two-dimensional distribution as in the precedent case, so that  we
can also apply Proposition~\ref{prop1B}.}
\end{remark}

\item $\alpha=0$, $\beta,\gamma\ne 0$, then we have that $e^1$, $\langle
e_3\rangle$ and $w$ determine $G_{1B}$ as the isotropy group.
Therefore, a $G_{1B}$-structure $P$ is given by an one-form $\omega$
without zeros, a one-dimensional distribution $L$ such
that $\omega|_{L}=0$, and a volume form $\Omega$.

\begin{proposition}
$P$ is integrable if and only if $\omega$ is closed.
\end{proposition}
{\bf Proof:}
If $d\omega=0$ then we have that the two-dimensional distribution
$D=\{\omega=0\}$ is involutive, and, by a similar argument as in Proposition
6.1, there are local coordinates $(y^1,y^2,y^3)$ such that
$\displaystyle{L=\langle\frac{\partial}{\partial y^3}\rangle}$,
$\displaystyle{D=\langle \frac{\partial}{\partial
y^2},\frac{\partial}{\partial y^1}\rangle}$ and
$\Omega=dy^1\wedge dy^2\wedge dy^3$.
Thus $\omega=\lambda dy^1$. Since $\omega$ is closed, we deduce that
$\lambda=\lambda(y^1)$ and, by defining
$x^1=\int\lambda(y^1)dy^1$,
$\displaystyle{x^2=\frac{y^2}{\lambda(y^1)}}$ and $x^3=y^3$,
we obtain adapted coordinates $(x^1,x^2,x^3)$. The converse is trivial.
\hfill $\Box$

\item $\beta=0$, $\alpha,\gamma\ne 0$,then $G_{1B}$ is
the isotropy group of the vector subspace $\langle e_1,e_3\rangle$,
the tangent covector $e_2$ in this subspace, and the tensor $w$.
Thus, a $G_{1B}$-structure  $P$ is given by a two-dimensional
distribution $D$, a tangent one-form $\omega|_{D}$ on $D$, and a volume
form $\Omega$.

\begin{proposition}
$P$ is integrable if and only if $D$ is involutive and $d\omega|_D=0$.
\end{proposition}

{\bf Proof:}
As in Proposition 6.6.
\hfill $\Box$

\item $\gamma=0$, $\alpha,\beta\ne 0$, then the group $G_{1B}$ is the
isotropy group of $e_3$, the subspace $\langle e_2,e_3\rangle$ and the
tensor $w$. Therefore, in this case, a $G_{1B}$-structure $P$ is given by
a two-dimensional distribution $D$, a vector field $X$ without zeros
belonging to $D$, and a volume form $\Omega$.

\begin{proposition}
$P$ is integrable if and only if $D$ is involutive and ${\cal L}_X\Omega=0$.
\end{proposition}

{\bf Proof:}
As in Proposition~\ref{prop1A}.
\hfill $\Box$
\end{itemize}

\bigskip

{\bf The group 1C}

\medskip

The group $G_{1C}$ is the isotropy group of the tensors $e_3,\,e^1$ and
$w$ on $\R^3$. Therefore a $G_{1C}$-structure $P$ is given by a one-form
$\omega$, a vector field $X$ such that $\omega(X)=0$, and a volume form
$\Omega$.

\begin{proposition}
$P$ is integrable if and only if $\omega$ is closed and ${\cal L}_X\Omega=0$.
\end{proposition}

{\bf Proof:} If $d\omega=0$ then the two-dimensional distribution
$D=\{\omega=0\}$ is involutive. Take a vector field $Y$ such that
$D=\langle X,Y\rangle$. Since $D$ is involutive, there are
functions $f$ and $g$, $g\ne 0$, such that $Y'=fX+gY$ verifies
$[X,Y']=0$. Thus, there exist local coordinates $(y^1,y^2,y^3)$
such that $\displaystyle{X=\frac{\partial}{\partial y^2}}$ and
$\displaystyle{Y'=\frac{\partial}{\partial y^3}}$. We deduce that
$\omega=\lambda dy^1$. Since $\omega$ is closed,
$\lambda=\lambda(y^1)$ and we can assume without loss of
generality that $\lambda=1$ (otherwise we only need to make a new
change of coordinates). Now $\Omega=b(y^1,y^2,y^3)dy^1\wedge
dy^2\wedge dy^3$, and by applying the argument in
Proposition~\ref{prop1A} we conclude. The converse is trivial.
\hfill $\Box$

\bigskip

{\bf The group 2A}

\medskip

The group $G_{2A}$ is the isotropy group of the subspaces $\langle
e_3\rangle$, $\langle e_2\rangle$ and of the tensor $w$. Thus, a
$G_{2A}$-structure $P$ is given by two transverse one-dimensional
distributions $L_1$ and $L_2$ and a volume form $\Omega$.

Assume that the two-dimensional distribution
$L_1\oplus L_2$ is involutive.
Given an adapted local basis $Y_1,Y_2,Y_3$, i.e., such that $Y_i\in
L_i$ for $i=1,2$ and $\Omega(Y_1,Y_2,Y_3)=1$, we define
\[
\tau(Y_{1}, Y_{2}, Y_{3}) =
Y_1(\alpha) + Y_2(h) + \alpha h - \alpha \beta,
\]
where $[Y_1,Y_2] = \alpha Y_1 + \beta Y_2$ and ${\cal
L}_{Y_1}\Omega=h\Omega$. Given another adapted local basis
$Y'_{1}=f_{1}Y_{1}$, $Y'_{2} = f_{2}Y_{2}$, and $\displaystyle{Y'_{3} =
\frac{1}{f_{1}f_{2}} Y_{3}}$, a direct computation shows that
\[
\tau(Y'_{1}, Y'_{2}, Y'_{3}) = f_{1} f_{2}
\tau(Y_{1}, Y_{2}, Y_{3}).
\]

\begin{proposition}\label{prop2A}
$P$ is integrable if and only if the two-dimensional distribution
$L_1\oplus L_2$ is involutive and $\tau (Y_{1},Y_{2},Y_{3})$
vanishes for an arbitrary adapted local frame $Y_{1}, Y_{2},
Y_{3}$.
\end{proposition}

{\bf Proof:}
One direction is obvious. For the other direction,
we only need to prove that given an adapted local frame $Y_1, Y_2,Y_3$
as above there exist local
functions $\lambda_1$ and $\lambda_2$ such that $[\lambda_1
Y_1,\lambda_2 Y_2]=0$ and ${\cal L}_{\lambda_1 Y_1}\Omega=0$. To prove
this we proceed as follows. If $[Y_1,Y_2]=\alpha Y_1+\beta Y_2$ and
${\cal L}_{Y_1}\Omega=h\Omega$, then
\[
Y_1(\lambda_2)+\beta\lambda_2=0,\hspace{5mm}Y_2(\lambda_1)-\alpha\lambda_1
=0,\hspace{5mm}Y_1(\lambda_1)+h\lambda_1=0.
\]
The first equation can always be integrated and the compatibility
condition for the last two equations is just
\begin{eqnarray*}
(-\alpha h+\beta\alpha)\lambda_1 & = & \alpha Y_1(\lambda_1)+\beta
Y_2(\lambda_1) \\
& = & [Y_1,Y_2]\lambda_1\\
& = & Y_1Y_2\lambda_1-Y_2Y_1\lambda_1 \\
& = & Y_1(\alpha \lambda_1)+Y_2(h\lambda_1)\\
& = & (Y_1\alpha+Y_2h)\lambda_1+\alpha
Y_1\lambda_1+hY_2\lambda_1\\
& = & (Y_1\alpha+Y_2 h)\lambda_1.
\end{eqnarray*}
To end the proof, we remark that if $[\lambda_1 Y_1,\lambda_2
Y_2]=0$ and ${\cal L}_{\lambda_1 Y_1} \Omega=0$, then there exist
local coordinates $(y^1,y^2,y^3)$ such that
$\displaystyle{\lambda_1 Y_1=\frac{\partial}{\partial y^1}}$,
$\displaystyle{\lambda_2 Y_2 = \frac{\partial}{\partial y^2}}$,
$\displaystyle{{\cal L}_{\frac{\partial}{\partial y^1}}\Omega=0}$,
$\displaystyle{{\cal L}_{\frac{\partial}{\partial y^2}}\Omega=0}$,
and now, after an appropriate change of coordinates, we conclude.
\hfill $\Box$

\bigskip

{\bf The group 2B}

\medskip

We can only consider the algebraic case, i.e., with $\alpha,\, \beta$
and $\gamma$ integer parameters.

\begin{itemize}
\item
If $\alpha\beta\gamma\ne 0$ then we have a tangent $H$-structure
defined by tangent tensors to the vector subspace $\langle
e_1,e_2\rangle$.
We can consider the following two subcases:
\begin{itemize}
\item If $\beta$ and $\gamma$ are both positive, we have
\[
\underbrace{e^1\otimes\cdots\otimes
e^1}_{\gamma}\otimes\underbrace{e_2\otimes\cdots\otimes e_2}_{\beta}\
\rm{and}\ \underbrace{e^2\otimes\cdots\otimes
e^2}_{\beta}\otimes\underbrace{e_1\otimes\cdots\otimes e_1}_{\gamma}
\]
as invariant tangent tensors of $G_{2B}$.

\item If $\beta>0$ and $\gamma<0$, then
\[
(\underbrace{e_1\otimes\cdots\otimes e_1}_{-\gamma})\odot
(\underbrace{e_2\otimes\cdots\otimes e_2}_{\beta})
\]
is an invariant tangent tensor.
\end{itemize}

Thus, the group $G_{2B}$ is the intersection of $G_{2A}$ with the
isotropy group of the above tangent tensors. Therefore, a
$G_{2B}$-structure $P$ is given by a one-dimensional distribution $L$,
a two-dimensional distribution $D$ with $L\subset D$, a tangent tensor
field on $D$ which is 0-deformable to the above tensors and a volume form.

>From Proposition~\ref{prop3.8}  and Corollaries~\ref{cor3.9} and
\ref{cor3.11}, we have the following.

\begin{proposition}\label{prop2B}
$\overline{P}$ is integrable if and only if $D$ is involutive and
the $H$-structures induced on the leaves of the foliation defined
by $D$ are all of them integrable (the last condition occurs if
and only if the unique torsion free $H$-connection on each leaf
has zero curvature). Assume that $\overline{P}$ is integrable with
adapted local coordinates $(x^1,x^2,x^3)$. If $P$ is integrable
then $\displaystyle{{\cal L}_{\frac{\partial}{\partial
x^2}}\Omega={\cal L}_{\frac{\partial}{\partial x^3}}\Omega=0}$.
\end{proposition}

{\bf Proof:} It only remains to prove the last assertion, which
follows from Remark~\ref{bk} and the computations made in the
proof of Proposition~\ref{prop1B}. \hfill $\Box$

\item If $\alpha\beta\gamma=0$ we consider the following subcases:

\begin{itemize}

\item $\alpha=0$, a $G_{2B}$-structure $P$ is given by two
one-dimensional distributions $L_1,L_2$, a $1$-form
$\omega$ such that $L_1\oplus L_2=\{\omega=0\}$ and a volume form
$\Omega$.

\begin{proposition}
$P$ is integrable if and only if the conditions in
Proposition~\ref{prop2A} hold and $d\omega=0$.
\end{proposition}

{\bf Proof:}
>From Proposition~\ref{prop2A} we conclude that there exist local coordinates
$x^1,x^2,x^3$ adapted to corresponding  $G_{2A}$-structure, this
implies that $\omega=\lambda dx^1$. Since $d\omega=0$, after an
appropriate change of coordinates we can assume that $\lambda=1$
and $\Omega = dx^1\wedge dx^2\wedge dx^3$. \hfill $\Box$

\item $\beta=0$, which is equivalent to case $\gamma=0$.
A $G_{2B}$-structure $P$ is given by a one-dimensional distribution
$L$, a vector field $X$ and a volume form $\Omega$.

\begin{proposition}\label{prop2BB}
$P$ is integrable if and only if
\[
{\cal L}_X L\subset L\hspace{3mm} \hbox{and} \hspace{3mm}
{\cal L}_X\Omega=0.
\]
\end{proposition}

\end{itemize}
\end{itemize}

\bigskip

{\bf The group 2C}

\medskip

The group $G_{2C}$ is the isotropy group of the vectors $e_2$, $e_3$ and
the tensor $w$. Thus a $G_{2C}$-structure $P$ is given by two vector fields
$X_1$, $X_2$ which are linearly independent,
and a volume form $\Omega$. We notice that $P$ can be alternatively
described by $X_1$, $X_2$ and a one-form $\omega$ such that
$\omega(X_i)=0$, $i=1,2$.
In fact, given $\Omega$ we put
$\omega=\iota_{X_1}\iota_{X_2}\Omega$, and given $\omega$, we define
$\Omega$ by $\Omega(X_1,X_2,Z)=\omega(Z)$ for all vector field $Z$.

\begin{proposition}
The following statements are equivalent:
\begin{enumerate}
\item $P$ is integrable
\item $[X_1,X_2]=0$ and ${\cal L}_{X_1}\Omega={\cal L}_{X_2}\Omega=0$.
\item $[X_1,X_2]=0$ and $d\omega=0$.
\end{enumerate}
\end{proposition}

\bigskip

{\bf The group 3A}

\medskip

The group $G_{3A}$ is the isotropy group of the subspaces $\langle
e_1,e_2\rangle$, $\langle e_2,e_3\rangle$ and the tensor $w$. Therefore
a $G_{3A}$-structure $P$ is given by two distributions of dimension two
$D_1$, $D_2$ and a volume form $\Omega$.

\begin{proposition}\label{prop3A}
$P$ is integrable if and only if $D_1$ and $D_2$ are both involutives.
\end{proposition}
{\bf Proof:}
Assume that $D_1$ and $D_2$ are both involutive. Then, by applying
Proposition~\ref{prop1A} to $L=D_1\cap D_2\subset D_1$ and $\Omega$, we
obtain local coordinates $(y^1,y^2,y^3)$ such that
$\displaystyle{D_1=\langle\frac{\partial}{\partial y^1},
\frac{\partial}{\partial y^2}\rangle}$ and $\displaystyle{D_1\cap D_2 =
\langle\frac{\partial}{\partial y^3}\rangle}$.
So, $\displaystyle{D_2=\langle\frac{\partial}{\partial y^3},
\frac{\partial}{\partial y^2} +
a\frac{\partial}{\partial y^1}\rangle}$, for some function $a$.
Since $D_2$ is involutive, we deduce that $\displaystyle{\frac{\partial
a}{\partial y^3}=0}$.
It is clear that
\[
\lambda \frac{\partial}{\partial y^1} + \mu\frac{\partial}{\partial
y^3},\, \frac{\partial}{\partial y^3},\, \frac{\partial}{\partial y^2}
+ a\frac{\partial}{\partial y^1}
\]
is an adapted local frame for all functions $\lambda$ and $\mu$,
$\lambda\ne 0$. Now, the equation
\[
\left[\lambda\frac{\partial}{\partial y^1}
+ \mu\frac{\partial}{\partial y^3},
\frac{\partial}{\partial y^2}
+ a\frac{\partial}{\partial y^1}\right]=0
\]
is equivalent to the following system of PDE's:
\begin{equation}\label{ece1}
\left(\frac{\partial}{\partial y^2} + a \frac{\partial}{\partial y^1}
\right)\lambda = -\frac{\partial a}{\partial y^1}\lambda \; , \;
\left(\frac{\partial}{\partial y^2}+a\frac{\partial}{\partial
y^1}\right)\mu = 0.
\end{equation}
In addition, we have that
$\displaystyle{\left[\frac{\partial}{\partial y^3},
\frac{\partial}{\partial y^2} + a \frac{\partial}{\partial
y^1}\right]=0}$ if and only if $\lambda$ and $\mu$ do not depend
on $y^3$. Thus~\mbox{(\ref{ece1})} have solution which are
independent of $y^3$ if and only if $a$ does not depend on $y^3$.
In such case, the other two brackets
\[
\left[
\frac{\partial}{\partial y^3},
\lambda\frac{\partial}{\partial y^1}+\mu\frac{\partial}{\partial y^3}
\right] \;\;{\rm and }\;\;
\left[
\frac{\partial}{\partial y^3},
\frac{\partial}{\partial y^2}+a\frac{\partial}{\partial y^1}
\right]
\]
also vanish and we conclude. The converse is trivial. \hfill
$\Box$

\bigskip

{\bf The group 3B}

\medskip

We only consider the algebraic case. This case is similar to the
case 2B but now the tensor fields are transverse to the
one-dimensional distribution $L$ instead of tangent to the
two-dimensional distribution $D$. Concerning the integrability we
obtain the following result.

\begin{proposition}
$\overline{P}$ is integrable if and only if it is projectable,
i.e. the coefficients of the transverse tensor field $T$ in a
system of coordinates $(x^1,x^2,x^3)$ adapted to the foliation
induced by $L$  do not depend on the tangent coordinate $x^3$.
Assume that $\overline{P}$ is integrable with adapted local
coordinates $(x^1,x^2,x^3)$. Then $P$ is integrable if and only if
$\displaystyle{{\cal L}_{\frac{\partial}{\partial
x^1}}\Omega={\cal L}_{\frac{\partial}{\partial x^2}}\Omega={\cal
L}_{\frac{\partial}{\partial x^3}}\Omega=0}$.
\end{proposition}

{\bf Proof:} The last assertion is a consequence of
Proposition~\ref{prop3.12}, Remark~\ref{bk} and the following
computation of $\overline{\mathfrak{g}}^{(1)}$: If
$\tau=(\tau_i)_{i=1}^3$ with $\tau_i=\left(\begin{array}{ccc}
\alpha
a_i & b_i & c_i\\ 0 & \beta a_i & 0\\
0 & 0 & \gamma a_i\end{array}\right)\in\overline{\mathfrak{g}}$,
then
\begin{center}
\begin{tabular}{c|c|c|c|}
~ & $\tau_{12}^k-\tau_{21}^k$ & $\tau_{13}^k-\tau_{31}^k$ &
$\tau_{23}^k-\tau_{32}^k$\\
\hline
 $k=1$ & $b_1-\alpha a_2$ & $c_1-\alpha a_3$ & $c_2-b_3$\\
\hline
 $k=2$ & $\beta a_1-0$ & $0-0$ & $0-\beta a_3$\\
\hline
 $k=3$ & $0-0$ & $\gamma a_1-0$ & $\gamma a_2-0$\\
\hline
\end{tabular}
\end{center}
Consequently, if $\tau\in\overline{\mathfrak{g}}^{(1)}$ then
$a_1=a_2=a_3=0$ and therefore
$\mathrm{tr}\,\tau_1=\mathrm{tr}\,\tau_2=\mathrm{tr}\,\tau_3=0$.
\hfill $\Box$

\bigskip

{\bf The group 3C}

\medskip

The group $G_{3C}$ is the isotropy group of the covectors $e^2$, $e^3$
and the tensor $w$. Thus a $G_{3C}$-structure $P$ is given by two
one-forms $\omega_1$ and $\omega_2$ which are linearly independent, and
a volume form $\Omega$. As in the case 2C we have an alternative
description of $P$, giving $\omega_1$,
$\omega_2$ and a vector field $X$ such that $\omega_1(X)=\omega_2(X)=0$.

\begin{proposition}
$P$ is integrable if and only if $\omega_1$ and $\omega_2$ are closed.
\end{proposition}
{\bf Proof:} From Proposition \ref{prop3A} we conclude that there
exist local coordinates $y^1,y^2,y^3$ such that
$\omega_i=\lambda_i dy^i$ for $i=1,2$ and $\Omega=dy^1\wedge
dy^2\wedge dy^3$. Since $d\omega_i=0$ then $\lambda_i$ only
depends on $y^i$, after the change of coordinates given by
$x^i=\int\lambda_i(y^i)dy^i$, $i=1,2$ and
$\displaystyle{x^3=\frac{y^3}{\lambda_1\lambda_2}}$ we conclude.
The converse is trivial. \hfill $\Box$

\bigskip

{\bf The group 4A}

\medskip

The group $G_{4A}$ is the isotropy group of the subspaces $\langle
e_1,e_2\rangle$, $\langle e_2\rangle$, $\langle e_3\rangle$ and the
tensor $w$. Therefore a $G_{4A}$-structure is given by one
two-dimensional distribution $D$, two one-dimensional distributions
$L_1$ and $L_2$ such that $L_1\subset D$ and $L_2\cap D=0$, and a
volume form $\Omega$.

\begin{proposition}
If $P$ is integrable then $D$ is involutive and the induced
$G_{2A}$-structure is integrable.
\end{proposition}

\begin{remark}{\rm
The converse does not hold. For instance take on $\R^3$ the
$G_{4A}$-structure given by $L_1=\langle Y_1\rangle$, $L_2=\langle
Y_2\rangle$ and $\Omega$ such that $\Omega(Y_1,Y_2,Y_3)=1$ where
\[
Y_1=\frac{\partial}{\partial y^1},\;\;
Y_2=\frac{\partial}{\partial y^2},\;\;
Y_3=y^2\frac{\partial}{\partial y^1}+\frac{\partial}{\partial
y^2}+\frac{\partial}{\partial y^3}.
\]
}
\end{remark}

\bigskip

{\bf The group 4B}

\medskip

The group $G_{4B}$ is the subgroup of $G_{2B}$ which leaves the
subspace $\langle e_1,e_2\rangle$ invariant. Thus, a $G_{4B}$-structure
consist of a $G_{2B}$-structure $(L,T,\Omega)$ and a two-dimensional
distribution $D'$ which is complementary of $L$.

A necessary condition for the integrability of a $G_{4B}$-structure is
given by Propositions~\ref{prop2B}-\ref{prop2BB}.

\bigskip

{\bf The group 4C}

\medskip

The group $G_{4C}$ is the isotropy group of the vectors $e_2$, $e_3$ and
the covector $e^1$. Therefore a $G_{4C}$-structure $P$ is given by two
vector fields which are linearly independent, and one-form $\omega$ such
that $\omega(X_1)=0$ and $\omega(X_2)=0$.

\begin{proposition}
$P$ is integrable if and only if $[X_1,X_2]=0$ and $d\omega=0$.
\end{proposition}

\bigskip

{\bf The group 5A}

\medskip

The group $G_{5A}$ is the isotropy group of the subspaces $\langle
e_1\rangle$, $\langle e_2\rangle$, $\langle e_3 \rangle$ and the tensor
$w$. Alternatively, $G_{5A}$ can be described as the isotropy group of
a diagonalizable endomorphism $f$ with three distinct eigenvalues, and
$w$. Then, a $G_{5A}$-structure $P$ is given by three one-dimensional
distributions $L_1$, $L_2$ and $L_3$, and a volume form $\Omega$.
Alternatively, $P$ can be described by a tensor field $h$ of type
$(1,1)$ which is 0-deformable to $f$, and the volume form $\Omega$.

\begin{remark}{\rm
The Lie algebra $\mbox{\fr g}$ is of finite type, indeed $\mbox{\fr
g}^{(1)}=0$, and then there is at most a free torsion
$G_{5A}$-connection $\nabla$. In fact, since the Nijenhuis tensor $N_{h}$ of
$h$ vanishes one can constructs local $G_{5A}$-connections without
torsion and, since the uniqueness, they coincide on the overlappings.
}
\end{remark}

\begin{proposition}
Let $P$ be a $G_{5A}$-structure and $\overline{P}$ an associated
$\overline{G}_{5A}$-structure. Then the following statements are
equivalent:
\begin{enumerate}
\item $\overline{P}$ is integrable.
\item The distributions $L_i\oplus L_j$ are both involutive.
\item $N_h=0$.
\end{enumerate}
Moreover, $P$ is integrable if and only if $N_h=0$ and the
$G_{5A}$-connection $\nabla$ has zero curvature.
\end{proposition}

\bigskip

{\bf The group 5B}

\medskip

We only consider the algebraic case, i.e., $\alpha,\beta$ and
$\gamma$ are integer parameters. Reordering if it is necessary, we can
assume that $\beta\ge 0$ and $\alpha\le 0$. Then, the tensor
\[
t = \underbrace{e^1\otimes\cdots\otimes e^1}_{\beta} \otimes
\underbrace{e_2\otimes\cdots\otimes e_2}_{-\alpha}
\]
is invariant by $G_{5B}$. In fact, $G_{5B}$ is the isotropy group of
the subspaces $\langle e_i\rangle$, $i=1,2,3$, and the tensors $t$ and
$w$. Therefore, a $G_{5B}$-structure $P$ is given by three
complementary one-dimensional distributions $L_i$, a tensor field $T$
which is 0-deformable to $t$, and a volume form $\Omega$.
According to the precedent section we obtain that if
$X_1,X_2,X_3$ is an adapted local frame to $P$ and we put
$[X_i,X_j]=\gamma_{ij}^kX_k$, then the integrability of $P$
implies the following conditions:
\begin{equation}\label{*1}
\gamma_{12}^3=0,\quad \gamma_{13}^2=0,\quad\gamma_{23}^1=0
\end{equation}
and
\begin{equation}\label{*2}
\gamma\gamma_{12}^1+\alpha\gamma_{23}^3=0, \quad
\gamma\gamma_{12}^2-\beta\gamma_{13}^3,
\quad\beta\gamma_{13}^1-\alpha\gamma_{23}^2=0.
\end{equation}

\begin{remark}{\rm
\begin{enumerate}
\item Equations (\ref{*1}) imply that the distributions
$L_i\oplus L_j$ are both involutive.
\item In this case $\partial_{\mbox{\fr g}}$ is injective, therefore if
Equations (\ref{*1}) and (\ref{*2}) hold, then there exist
a unique torsionless $G_{5B}$-connection $\nabla$.
\end{enumerate}
}
\end{remark}

\begin{proposition}
$P$ is integrable if and only if $L_i\oplus L_j$ is involutive,
Equations (\ref{*2}) hold for any local adapted frame to $P$, and the
curvature of $\nabla$ vanishes.
\end{proposition}

\begin{remark}{\rm
The case $\alpha=0$ corresponds to a subgroup of $G_{19}$ which
is conjugated with the special Lorentz group as we will see later. In
this case, a $G_{5B}$-structure is given by a Lorentzian metric $g$,
its associated volume form, and a vector field $X$ such that
$g(X,X)=1$. Then $P$ is integrable if and only if the Ricci
tensor of $g$ is zero and $\nabla X=0$, where $\nabla$ is the
Levi-Civita connection of $g$.
}
\end{remark}

\bigskip

{\bf The group 5C}

\medskip

The group $G_{5C}$ is the trivial group. Then a $G_{5C}$-structure $P$
is just a linear parallelism $X_1,X_2,X_3$ on $B$.

\begin{proposition}
$P$ is integrable if and only if
$[X_1,X_2]=[X_1,X_3]=[X_2,X_3]=0$, or, equivalently, if and only if the
flat connection defined by the parallelism is symmetric.
\end{proposition}

\bigskip

{\bf The group 6A}

\medskip

The group $G_{6A}$ is the isotropy group of the vector subspace
$\langle e^1\otimes e^1+e^2\otimes e^2\rangle\subset T^2_0\R^3$ and
$w$. Thus, giving a $G_{6A}$-structure $P$ is equivalent to giving the
projectivization of a symmetric covariant tensor field of order 2 and
constant rank 2, and a volume form. But this is
equivalent to give a one-dimensional distribution $L$ and a transverse
almost complex structure $J$.
We denote by $\overline{P}$ the $\overline{G}_{6A}$-structure obtained
from $P$ without considering the volume form $\Omega$. Using the fact
that all $GL(1,\C)$-structure is integrable we obtain the following.

\begin{proposition}
$\overline{P}$ is integrable if and only if the transverse almost
complex structure $J$ is projectable, i.e., the transverse tensor
field $J$ does not depend on the tangent coordinates adapted to the
foliation defined by $L$.
\end{proposition}

\bigskip

{\bf The group 6B}

\medskip

The only algebraic subgroup of type 6B is obtained with $\alpha=0$.
In this case, a $\overline{G}_{6B}$-structure is given by a
one-dimensional distribution $L$ and a Riemannian metric $g$
which is transverse to $L$. If in addition, we give a volume form
$\Omega$ we obtain the corresponding $G_{6B}$-structure $P$.

\begin{proposition}
$\overline{P}$ is integrable if and only if $g$ is projectable, i.e.,
$g$ is a transverse bundle like metric to the foliation defined by $L$.
\end{proposition}

\bigskip

{\bf The group 7A}

\medskip

The group $G_{7A}$ is the isotropy group of the subspace $\langle
e_1\otimes e_1+e_2\otimes e_2\rangle\subset T_2^0\R^3$ and the tensor
$w$. Thus, a $G_{7A}$-structure $P$ is given by the projectivization of
a symmetric two-contravariant tensor field $g$ of rank 2, and a volume
form $\Omega$. An alternative description of $\overline{P}$ consists of
a two-dimensional distribution $D$ and a tangent almost complex structure
on $D$. From Proposition~\ref{prop3.8} we obtain the following result.

\begin{proposition}
$\overline{P}$ is integrable if and only if $D$ is involutive.
\end{proposition}

\bigskip

{\bf The group 7B}

\medskip

The only algebraic subgroup of type 7B is obtained by putting
$\alpha=0$. This group $G_{7B}$ is the isotropy group of the subspace
$\langle e_1,e_2\rangle$, the tangent metric tensor $e^1\otimes
e^1+e^2\otimes e^2$ and the tensor $w$. Then, a $G_{7B}$-structure $P$
is given by a two-dimensional distribution $D$ with a tangent metric
$g$ and a volume form $\Omega$.

\begin{proposition}\label{prop6.29}
$\overline{P}$ is integrable if and only if $D$ is involutive and
the scalar curvature of the metric defined on each leaf of the
induced foliation vanishes. Assume that $\overline{P}$ is
integrable with adapted local coordinates $(x^1,x^2,x^3)$. If $P$
is integrable then $\displaystyle{{\cal
L}_{\frac{\partial}{\partial x^1}}\Omega={\cal
L}_{\frac{\partial}{\partial x^2}}\Omega=0}$.
\end{proposition}

{\bf Proof:} The last assertion is a consequence of
Remark~\ref{bk} and the following calculation of
$\overline{\mathfrak{g}}^{(1)}$: If $\tau=(\tau_i)_{i=1}^3$ with
$\tau_i=\left(\begin{array}{ccc} 0 & a_i & b_i\\  -a_i & 0 & c_i\\
0 & 0 & \beta_i
\end{array}\right)\in\overline{\mathfrak{g}}$,
then
\begin{center}
\begin{tabular}{c|c|c|c|}
~ & $\tau_{12}^k-\tau_{21}^k$ & $\tau_{13}^k-\tau_{31}^k$ &
$\tau_{23}^k-\tau_{32}^k$\\
\hline
 $k=1$ & $a_1-0$ & $b_1-0$ & $b_2-a_3$\\
\hline
 $k=2$ & $0+a_2$ & $c_1+a_3$ & $c_2-0$\\
\hline
 $k=3$ & $0-0$ & $\beta_1-0$ & $\beta_2-0$\\
\hline
\end{tabular}
\end{center}
Consequently, if $\tau\in\overline{\mathfrak{g}}^{(1)}$ then
$\beta_1=\beta_2=0$ and therefore
$\mathrm{tr}\,\tau_1=\mathrm{tr}\,\tau_2=0$. \hfill $\Box$

\bigskip

{\bf The group 8A}

\medskip

The group $G_{8A}$ is the isotropy group of the endomorphism
$f=e^2\otimes e_1-e^1\otimes e_2$ and the tensor $w$. Thus, a
$\overline{G}_{8A}$-structure is given by a tensor field $h$ of type
(1,1) which is 0-deformable to $f$. We notice that $h^3+h=0$. Theses
structures are called $f$-structures in the literature \cite{Yano}.
If, in addition, we give a volume form $\Omega$ we obtain
the corresponding $G_{8A}$-structure $P$.

\begin{proposition}
$\overline{P}$ is integrable if and only if $N_h=0$. In this case,
there exist adapted local coordinates $x^1,x^2,x^3$ such that
$\Omega = b(x^1,x^2,x^3) \, dx^1 \wedge dx^2 \wedge dx^3$, for
some function $b$. Therefore, $P$ is integrable if and only if the
following conditions hold:
\begin{equation}\label{compatib}
 \frac{\partial^2\log b}{\partial
x^1\partial x^3}=0, \quad \frac{\partial^2\log b}{\partial
x^2\partial x^3}=0, \quad \frac{\partial^2 \log b}{\partial
(x^1)^2} + \frac{\partial^2\log b}{\partial (x^2)^2}=0.
\end{equation}
\end{proposition}

{\bf Proof:} The last assertion follows by applying the techniques
described at the end of Section~3. By some calculations similar to
the ones made in the proof of Proposition~\ref{prop6.29}, we
deduce that
\[\overline{\mathfrak{g}}^{(1)}=\left\{\tau_1=\left(\begin{array}{ccc}
\alpha_1 & \alpha_2 & 0\\ -\alpha_2 & \alpha_1 & 0\\ 0 & 0 &0
\end{array}\right),\ \tau_2=\left(\begin{array}{ccc}
\alpha_2 & -\alpha_1 & 0\\ \alpha_1 & \alpha_2 & 0\\ 0 & 0 &0
\end{array}\right),\ \tau_3=\left(\begin{array}{ccc}
0 & 0 & 0\\ 0 & 0 & 0\\ 0 & 0 &\beta_3
\end{array}\right) \right\}.\]
In order to obtain the integrability of $P$ we need to construct a
tensor field $\tau:U\to\overline{\mathfrak{g}}^{(1)}$ such that
$\frac{\partial \tau_j}{\partial x^i}-\frac{\partial
\tau_i}{\partial x^j}+[\tau_i,\tau_j]=0$ and
$\mathrm{tr}\,\tau_i=\frac{\partial\log b}{\partial x^i}$, where
$\Omega=b(x^1,x^2,x^3)dx^1\wedge dx^2\wedge dx^3$. Taking into
account that $[\tau_i,\tau_j]=0$, $\frac{\partial \log b
}{\partial x^i}=\mathrm{tr}\,\tau_i=2\alpha_i$ if $i=1,2$ and
$\frac{\partial\log b}{\partial x^3}=\mathrm{tr}\,\tau_3=\beta_3$,
the compatibility relations of the  resulting system of PDE's can
be expressed as (\ref{compatib}).\hfill $\Box$

\bigskip

{\bf The group 8B}

\medskip

The only algebraic subgroup of type 8B is obtained when $\alpha=0$. In
this case $G_{8B}$ is the subgroup of $SO(3)$ that leaves invariant the
vector $e_3$. Thus, a $G_{8B}$-structure $P$ is given by a Riemannian
metric $g$, the Riemannian volume form $\Omega$, and a vector field
$X$ without zeros.

\begin{proposition}
$P$ is integrable if and only if the Ricci tensor of $g$ vanishes and
$\nabla X=0$ where $\nabla$ is the Levi-Civita connection of $g$.
\end{proposition}

\bigskip

{\bf The group 9}

\medskip

The group $G_9$ is the special linear group $Sl(3,\R)$. Then a
$G_9$-structure is given by a volume form $\Omega$.
If we look in the proof of Proposition~\ref{prop1A},
we deduce the following result.

\begin{proposition}
Every $Sl(3,\R)$-structure is integrable.
\end{proposition}

\bigskip

{\bf The group 10}

\medskip

The group $G_{10}$ is the isotropy group of the subspace $\langle e_1,
e_2\rangle$ of $\R^3$ and the tensor $w$.
Then a $G_{10}$-structure $P$ is given by a two-dimensional
distribution $D$ and a volume form $\Omega$.

\begin{proposition}
$P$ is integrable if and only if the distribution $D$ is
involutive.
\end{proposition}

\bigskip

{\bf The group 11}

\medskip

The group $G_{11}$ is the isotropy group of the subspace $\langle
e_3\rangle$ of $\R^3$ and the tensor $w$. Thus, a $G_{11}$-structure
$P$ is given by a one-dimensional distribution $L$ and a volume form
$\Omega$.

\begin{proposition}
$P$ is always integrable.
\end{proposition}

\bigskip

{\bf The group 12}

\medskip

The group $G_{12}$ is the isotropy group of the covector  $e^3$ and
the tensor $w$. Therefore, a $G_{12}$-structure $P$ is given by a
one-form $\omega$ without zeros and a volume form
$\Omega$.

\begin{proposition}
$P$ is integrable if and only if $d\omega=0$.
\end{proposition}

\bigskip

{\bf The group 13}

\medskip

The group $G_{13}$ is the isotropy group of the vector $e_3$ and the
tensor $w$. Then, a $G_{13}$-structure $P$ is given by a vector field
$X$ without zeros and a volume form $\Omega$.

\begin{proposition}
$P$ is integrable if and only if ${\cal L}_{X}\Omega=0$.
\end{proposition}

\bigskip

{\bf The group 14}

\medskip

The group $G_{14}$ is the isotropy group of the subspaces $\langle
e_1,e_2\rangle, \langle e_3\rangle$ of  $\R^3$ and the tensor $w$.
Alternatively, it is the isotropy group of the endomorphism
$e^1\otimes e_1$ and the tensor $w$. Thus, a $G_{14}$-structure
$P$ is given by a two-dimensional distribution $D$, a
complementary one-dimensional distribution  $L$, and a volume form
$\Omega$. Alternatively, it is given by a tensor field $h$ of type
$(1,1)$ such that $h^3-h^2=0$, together with the form $\Omega$.
Let $\overline{P}$ be the $\overline{G}_{14}$-estructure obtained
from $P$ without considering the volume form $\Omega$.

\begin{proposition}
$\overline{P}$ is integrable if and only if $N_h=0$.
\end{proposition}

\bigskip

{\bf The group 15}

\medskip

The group $G_{15}$ is the isotropy group of the tensors $u=e^1\wedge
e^2$ and $v=e^3$. We notice that $w=u\wedge v$. Therefore a
$G_{15}$-structure $P$ is given by a $2$-form $\eta$ and a $1$-form
$\omega$ such that $\eta\wedge\omega \not= 0$ everywhere.
Theses structures are called almost-cosymplectic \cite{Blair}.

\begin{proposition}
$P$ is integrable if and only if $d\eta=0$ and $d\omega=0$.
\end{proposition}
{\bf Proof:}
It is just the statement of Darboux's Theorem (see \cite{Blair}).
\hfill $\Box$

\bigskip

{\bf The group 16}

\medskip

The group $G_{16}$ is the special orthogonal group $SO(3)$. Then a
$G_{16}$-structure $P$ is given by a Riemannian metric $g$ on $B$ and
the Riemannian volume form $\Omega$.

\begin{proposition}
$P$ is integrable if and only if the curvature of the Levi-Civita
connection of $g$ vanishes. This condition is equivalent to the
vanishing of the Ricci tensor of $g$.
\end{proposition}

\bigskip

{\bf The group 19}

\medskip

The group $G_{19}$ is the isotropy group of the tensors
\[
s = e^2\otimes e^3 + e^3\otimes e^2 - e^1\otimes e^1
\]
and $w$. We notice that $s$ is in the same orbit that $-e^1\otimes
e^1-e^2\otimes e^2+e^3\otimes e^3$, which is the model for the
Lorentzian metrics, thus the group $G_{19}$ is conjugated with the
special group $SO(2,1)$. Since we are discussing subgroups modulo
conjugation we can consider a $G_{19}$-structure as given by a Lorentz
metric $g$ and its volume form.

\begin{proposition}
As in the Riemannian case, $P$ is integrable if and only if the Ricci
tensor of $g$ is zero.
\end{proposition}

\bigskip

{\bf The group 20}

\medskip

The group $G_{20}$ is the isotropy group of the subspace of $T_0^4\R^3$
generated by
\[
t = (e^2 \otimes e^3 + e^3 \otimes e^2 - e^1 \otimes e^1) \otimes e^3
\otimes e^3 - e^3 \otimes e^3 \otimes (e^2 \otimes e^3 + e^3 \otimes
e^2 - e^1 \otimes e^1)
\]
and the tensor $w$. Therefore a $G_{20}$-structure $P$ is given by
the projectivization of a tensor field $T$ which is 0-deformable
to $t$, and a volume form $\Omega$. We can use the conditions in
Section 3, however we did not find any nice geometrical
interpretation.

\bigskip

{\bf The group 23}

\medskip

The group $G_{23}$ is the subgroup of $G_{19}$ which leaves the
subspace $\langle e_2\rangle$ invariant. Via conjugation of
$G_{19}$ with the special group $SO(2,1)$, $e_2$ becomes a
vector in the light cone. Thus, a $G_{23}$-structure $P$ is
given by a Lorentzian metric $g$, its volume form $\Omega$ and a
one-dimensional distribution $L$ contained in the light cone of $g$.

\begin{proposition}
$P$ is integrable if and only if the Ricci tensor of $g$ vanishes and
$\nabla L\subset L$, where $\nabla$ is the generalized Levi-Civita
connection of $g$.
\end{proposition}

\bigskip

{\bf The group 24}

\medskip

The group $G_{24}$ is conjugated with the isotropy group of the tensor
$w$ and the endomorphism $u=e^2\otimes e_1+e^3\otimes e_2$, with
minimum polynomial $u^3=0$. Then, a $\overline{G}_{24}$-structure
$\overline{P}$ is given by a tensor field $h$ which is 0-deformable to
$u$. If, in addition, we give a volume form $\Omega$ then we obtain the
corresponding $G_{24}$-structure $P$.

\begin{proposition}
$\overline{P}$ is integrable if and only if $N_h=0$.
Assuming that $\overline{P}$ is integrable, then we deduce that $P$ is
integrable if and only if for an arbitrary system of
coordinates $x^1,x^2,x^3$ adapted to $\overline{P}$ the following
conditions hold:
\begin{equation}\label{compatib2}
{\cal L}_{\frac{\partial}{\partial x^1}}\Omega=0, \quad {\cal
L}_{\frac{\partial}{\partial x^2}}\Omega=0.
\end{equation}
\end{proposition}

{\bf Proof:} Using the conjugation between $\overline{G}_{24}$ and
$\overline{G}=G_u$ we can characterize the integrability of $P$ in
terms of $\overline{\mathfrak{g}}=\left\{\left(\begin{array}{ccc}
a & b & c\\ 0 & a & b\\ 0 & 0 & a
\end{array}\right)\ |\ a,b,c\in\R\right\}$. First of all, we
compute $\overline{\mathfrak{g}}^{(1)}$ by means of the following
table
\begin{center}
\begin{tabular}{c|c|c|c|}
~ & $\tau_{12}^k-\tau_{21}^k$ & $\tau_{13}^k-\tau_{31}^k$ &
$\tau_{23}^k-\tau_{32}^k$\\
\hline
 $k=1$ & $b_1-a_2$ & $c_1-a_3$ & $c_2-b_3$\\
\hline
 $k=2$ & $a_1-0$ & $b_1-0$ & $b_2-a_3$\\
\hline
 $k=3$ & $0-0$ & $a_1-0$ & $a_2-0$\\
\hline
\end{tabular}
\end{center}
obtaining that
\[\overline{\mathfrak{g}}^{(1)}=\left\{\tau_1=\left(\begin{array}{ccc}
0 & 0 & a_3\\ 0 & 0 & 0\\ 0 & 0 & 0
\end{array}\right),\ \tau_2=\left(\begin{array}{ccc}
0 & a_3 & b_3\\ 0 & 0 & a_3\\ 0 & 0 & 0
\end{array}\right),\ \tau_1=\left(\begin{array}{ccc}
a_3 & b_3 & c_3\\ 0 & a_3 & b_3\\ 0 & 0 & a_3
\end{array}\right)\right\}.\]
The integrability of $P$ is equivalent to the existence of a
tensor field $\tau:U\to\overline{\mathfrak{g}}^{(1)}$ verifying
the system of PDE's $\frac{\partial \tau_j}{\partial
x^i}-\frac{\partial \tau_i}{\partial x^j}+[\tau_i,\tau_j]=0$,
jointly to the equations $\mathrm{tr}\,\tau_i=\frac{\partial\log
b}{\partial x^i}$, being $\Omega=b(x^1,x^2,x^3)dx^1\wedge
dx^2\wedge dx^3$. Taking into account that $[\tau_i,\tau_j]=0$, it
is easy to check that the precedent system of PDE's has always
solution provides that $a_3=\frac{1}{3}\frac{\partial\log
b}{\partial x^3}$ and the function $b(x^1,x^2,x^3)$ does not
depend on $x^1$ nor $x^2$, what is equivalent to
(\ref{compatib2}). \hfill $\Box$

\bigskip

{\bf The group 26}

\medskip

The group $G_{26}$ is another subgroup of $G_{19}$. Under the same
identifications as in that case, and by similar computations
we conclude that giving a $G_{26}$-structure $P$ is the same that
giving a Lorentzian metric $g$, its volume form $\Omega$, and a
isotropic vector field $X$ without zeros (i.e., such that $g(X,X)=0$).

\begin{proposition}
$P$ is integrable if and only if the Ricci tensor of $g$ vanishes and
$\nabla X=0$, where $\nabla$ is the generalized Levi-Civita connection
of $g$.
\end{proposition}

\section{Chevalley's Theorem}

As we have shown, a $G$-structure defined by a 0-deformable tensor
field is a reduction of the frame bundle ${\cal F}M$ to an
algebraic subgroup $G$ of $\GL{n}$. We can ask for a converse:
given a $G$-structure with $G$ an algebraic subgroup of $\GL{n}$,
there exist a 0-deformable tensor field defining it?

An approach to this question is given by a Theorem of Chevalley
whose proof we sketch below. Before stating it, let us introduce
some notions. Let $V$ be a vector space. A {\sf construction} over
$V$ is a vector space obtained from $V$ by iterating the
operations $*,\oplus,\otimes,S^m$ and $\bigwedge^m$. If $g\in
GL(V)$ then $g$ acts in a natural way on each construction over
$V$ following the rules:
\begin{itemize}
\item[(i)] if $v\in V$ then $\rho(g)(v)=gv$ is the standard
action,
\item[(ii)] if $\omega\in V^*$ then $\rho(g)(\omega):v\mapsto
\omega(g^{-1}v)$,
\item[(iii)] $\rho(g)(a\otimes b)=\rho(g)(a)\otimes\rho(g)(b)$,
\item[(iv)] $\rho(g)(a\oplus b)=\rho(g)(a)\oplus\rho(g)(b)$.
\end{itemize}
In other words, $\rho$ defines the (faithful) tensorial
representation $\rho:GL(V)\to GL(T_V)$, where
$T_V=\bigoplus_{r,s\ge 0}V^{r,s}$ (with $V^{r,s}=V^{\otimes
r}\otimes (V^*)^{\otimes s}$) is the whole tensor algebra over
$V$.

\begin{remark}{\rm
Any construction $W$ over $V$ is a direct sum of finite
dimensional subspaces of the whole tensor algebra $T_V$.
Therefore, if $M$ is a manifold of dimension $n=\dim V$, the
associated vector bundle $({\cal F}M\times W)/\GL{n}$ is a direct
sum of tensor bundles over $M$. On the other hand, if $u\in W$ and
$G=\{g\in\GL{n}\,|\, \rho(g)u=u\}$ (respectively,
$G=\{g\in\GL{n}\,|\, \exists\lambda\neq 0,\ \rho(g)u=\lambda u\}$)
then, according to Theorem~\ref{Teorema}, any $G$-structure is
given by some $W$-tensor (respectively $\mathbb{P}W$-tensor) which
is $0$-deformable to $u\in W$ (respectively,
$[u]\in\mathbb{P}W$).}
\end{remark}

It is easy to prove that if $\{V_i\}_{i\in I}$ is a family of
constructions, and for every $i\in I$, $\{W_{ij}\}_{j\in J_i}$ is
a family of subspaces of $V_i$, then
\[
H = \{ g\in GL(V)\,|\, g(W_{ij}) \subset W_{ij}, \forall i\in I,
\, \forall j \in J_i\}
\]
is an algebraic subgroup of $GL(V)$. The following result is a
weak converse.

\begin{theorem}[Chevalley](see \cite{Ber,Ramis})
Let $H$ be an algebraic subgroup of $GL(V)$. Then there exist a
finite-dimensional vector space $W$, a faithful representation
$\alpha:GL(V)\to GL(W)$ and a subspace $W_0\subset W$ such that
\[H=\{g\in GL(V)\,|\, \alpha(g)(W_0)\subset W_0\}.\]
In fact, $W_0$ and $W$ can be chosen such that $\dim W_0=1$.
%
\end{theorem}

{\bf Sketch of the proof (from \cite{Ber,Borel}):} For the sake of
shortness we will write $G$ instead of $GL(V)$ and we denote by
$K[G]$ the ring of regular functions over $G$ with values in the
ground field $K$ ($K=\R$ in what follows), i.e.
$K[G]=K[x_{ij},\,1\le i,j\le n][\Delta^{-1}]$ with
$\Delta=\det(x_{ij})$ and $n=\dim V$.

Since $G$ acts on itself by left translations, it induces a linear
action on $K[G]$ by pull-back. More precisely, if $g\in G$ and
$f\in K[G]$ then $g_*f:=(g^{-1})^*f:(x_{ij})\mapsto
f(g^{-1}\cdot(x_{ij}))$. The map $g\mapsto g_*$ defines a faithful
representation $\alpha:G\to GL(K[G])$. The space $W$ stated in the
theorem will be a suitable subspace of $K[G]$ with the restricted
representation $\alpha:G\to GL(W)$.

Since $H\subset G$ is algebraic, we can consider the ideal
$I(H)\subset K[G]$ defining $H$, i.e. $I(H)=\{f\in K[G]\,|\,
f(h)=0 \textrm{ for all }h\in H\}$. Using that $K[G]$ is
noetherian, we can find a finite dimensional subspace $W\subset
K[G]$ containing a system of generators of $I(H)$. Then, we define
$W_0$ as $W\cap I(H)$ and we have
\[H=\{g\in G\,|\, g^{-1}\cdot H\subset H\}=\{g\in G\,|\,\alpha(g)(I(H))\subset I(H)\}=\{g\in G\,|\, \alpha(g)(W_0)\subset W_0\}.\]
Finally, if $k=\dim W_0>1$ then we can take
$\alpha'=\bigwedge^k\alpha$, $W'=\bigwedge^k W$ and
$W_0'=\bigwedge^k W_0$ which is one-dimensional. \hfill $\Box$

\begin{remark}\label{tensor}
{\rm Denoting by $G=GL(V)$ and by $E=End(V)=V\otimes V^*$, notice
that $K[G]$ contains $K[x_{ij}]=K[E]\cong S(E^*)=\bigoplus_{s\ge
0}(E^*)^{\otimes s}\hookrightarrow T_V$. On the other hand, from
the proof we can choose $W\subset K[E]$. Thus, we can take
$W\subset T_V$, but the representation $\alpha:GL(V)\to GL(W)$ is
not the restriction of the tensorial representation $\rho:GL(V)\to
GL(T_V)$ considered above.}
\end{remark}

Although in some references \cite{Che,Oni} it seems that
Chevalley's theorem holds for $\alpha=\rho$, we have not found a
proof of this. Therefore, we prefer do not use this stronger
version and, consequently, we show, in the same spirit, a slightly
different result:
\begin{theorem}\label{thm7.4}
Let $H$ be an algebraic subgroup of $GL(V)$. Then there exists a
finite dimensional subspace $W\subset T_V$ such that the
normalizer ${\cal N}(H)=\{g\in GL(V)\,|\, g^{-1}Hg=H\}$ of $H$
verifies:
\[{\cal N}(H)=\{g\in GL(V)\,|\, \rho(g)(W)\subset W\},\]
where $\rho:GL(V)\to GL(T_V)$ denotes the tensorial
representation. Moreover, we can take $W$ one-dimensional.
\end{theorem}

{\bf Proof:} Consider the adjoint action of $G$ onto itself and
the induced representation $\rho':G\to GL(K[G])$ given by
$\rho'(g)(f):(x_{ij})\mapsto f(g^{-1}\cdot(x_{ij})\cdot g)$ for
any $f\in K[G]$. A straightforward calculation shows that, under
the identifications made in Remark~\ref{tensor}, $\rho'$ coincides
with the tensorial representation $\rho:G\to GL(T_V)$.

Now, we proceed as above by considering the ideal $I(H)$ of $K[G]$
defining the algebraic subgroup $H\subset G$. Again, since  $K[G]$
is noetherian there is a finite dimensional subspace $W$ of
$K[E]\subset K[G]$ which contains a system of generators of
$I(H)$. Identifying $W$ as a subspace of $T_V$, we have
\begin{eqnarray*}
&&\{g\in G\,|\,\rho(g)(W)\subset W\}=\{g\in G\,|\,\rho'(g)(I(H))\subset I(H)\}\\
&=&\{g\in G\,|\,\rho'(g)(f)\in I(H),\forall f\in I(H)\}=\{g\in
G\,|\, \rho'(g)(f)(h)=0,\forall f\in I(H),\forall h\in H\}\\
&=&\{g\in G\,|\, f(g^{-1}hg)=0,\forall f\in I(H),\forall h\in
H\}=\{g\in G\,|\, g^{-1}hg\in H,\forall h\in H\}\\
&=&{\cal N}(H).
\end{eqnarray*}
which proves the result. The last assertion follows taking a
suitable exterior power of $W$. \hfill $\Box$

\begin{theorem}
Let $G$ be an algebraic subgroup of $\GL{n}$. Then, any ${\cal
N}(G)$-structure can be given by the projectivization of an
inhomogeneous tensor field $t=\sum\limits_{i=1}^mt_i$, where each
$t_i$ is a $0$-deformable tensor field on $M$ of type $(r_i,s_i)$.
\end{theorem}

\begin{remark}
{\rm In \cite{EdL2000} the authors explore the possibility of the
existence of non-uniform materials which could enjoy, however,
some kind of homogeneity. They introduce the notion of
unisymmetric materials as follows.

\begin{definition}
A material body is said to be {\sf unisymmetric} if the material
symmetry groups of its points in one (and, therefore, in every)
configuration are pairwise conjugate.
\end{definition}

Functionally graded materials (FGM for short), important for their
industrial applications, are of this type.

Let $X_1$ and $X_2$ be two points of a unisymmetric body $B$, and
let $A : T_{X_1}B \longrightarrow T_{X_2}B$ be a symmetry
isomorphism such that $G_2 = AG_1 A^{-1}$, where $G_1$ and $G_2$
are the material symmetry groups at $X_1$ and $X_2$, respectively.
Then, the family ${\cal A}_{12}$ of all possible symmetry
isomorphisms between both points is
\begin{equation}\label{normal}
{\cal A}_{12} = A {\cal N}(G_1),
\end{equation}
where ${\cal N}(G_1)$ is the normalizer of $G_1$ in $\GL{n}$. If
we proceed as in Section 5 for uniform materials, and choose a
point $X_0$ and a particular linear frame $Z_0$ at $X_0$, we can
transport the material symmetry group $G(X_0)$ at $X_0$ to $\R^n$
and obtain a subgroup $G$ of $\GL{n}$. Using (\ref{normal}) we
deduce that possible admissible references at $X_0$ is just $Z_0
{\cal N}(G)$, where ${\cal N}(G)$ is the normalizer of $G$. This
means that the geometric structure associated to a unisymmetric
body is a ${\cal N}(G)$-structure. Accordingly to
Theorem~\ref{thm7.4} this implies that, if $G$ is an algebraic
subgroup of $\GL{n}$ then the geometric ${\cal N}(G)$-structure is
defined by the projectivization of a tensor. This fact probably
deserves a more careful analysis to be done elsewhere.
 }
\end{remark}

\section*{Acknowledgments}

This work has been supported through grants DGICYT (Spain),
Projects PB94-0106, PB97-1257, BFM2001-2272 and NATO CRG 950833.
D.M. acknowledges the support of the Consejo Superior de
Investigaciones Cient{\'\i}ficas through a grant. We acknowledge
to Profs. Marcel Nicolau and Agust{\'\i} Revent\'os for their help
during the preparation of this paper. M. de L. wish to acknowledge
the warm hospitality of the Departament de Matem\`atiques
(Universitat Aut\`onoma de Barcelona).

\newpage

\section*{A. Connected Lie subgroups of the $Sl(3,\R)$ }

\footnotesize

\begin{center}

\vspace{2mm}

\begin{tabular}{|cllc|}
\hline
Number&Lie subgroups with det$\,G=1$&Lie Subalgebras with tr$\,H=0$&Dimension\\
\hline
& & & \\
1&$\left[\begin{array}{ccc}
e^{\alpha a}&0&0\\
b&e^{\beta a }& 0\\
c&d&e^{\gamma a}
\end{array}\right]$&$\left[\begin{array}{ccc}
\alpha a& 0& 0\\
b &\beta a & 0\\
c&d&\gamma a
\end{array}\right]$&$
5,\,4,\,3$\\  & & & \\
2&$\left[\begin{array}{ccc}
e^{\alpha a}&0&0\\
b&e^{\beta a }& 0\\
c&0&e^{\gamma a}
\end{array}\right]$&$\left[\begin{array}{ccc}
\alpha a& 0& 0\\
b &\beta a & 0\\
c&0&\gamma a
\end{array}\right]$&$
4,\,3,\,2$\\
& & & \\
3&$\left[\begin{array}{ccc}
e^{\alpha a}&b&c\\
0&e^{\beta a }& 0\\
0&0&e^{\gamma a}
\end{array}\right]$&$\left[\begin{array}{ccc}
\alpha a& b&c\\
0 &\beta a & 0\\
0&0&\gamma a
\end{array}\right]$&$
4,\,3,\,2$\\
  & & & \\
4&$\left[\begin{array}{ccc}
e^{\alpha a}&0&0\\
b&e^{\beta a }& 0\\
0&0&e^{\gamma a}
\end{array}\right]$&$\left[\begin{array}{ccc}
\alpha a& 0& 0\\
b &\beta a & 0\\
0&0&\gamma a
\end{array}\right]$&$
3,\,2,\,1$\\  & & & \\
5&$\left[\begin{array}{ccc}
e^{\alpha a}&0&0\\
0&e^{\beta a }& 0\\
0&0&e^{\gamma a}
\end{array}\right]$&$\left[\begin{array}{ccc}
\alpha a& 0& 0\\
0 &\beta a & 0\\
0&0&\gamma a
\end{array}\right]$&$
2,\,1,\,0$\\ & & & \\
6&$\left[\begin{array}{ccc}
e^{\alpha a}\cos(a)&e^{\alpha a}\sin(a)&0\\
-e^{\alpha a}\sin(a)&e^{\alpha a}\cos(a)&0\\
b&c&e^{\beta a}
\end{array}\right]$&$\left[\begin{array}{ccc}
\alpha a &a&0\\
-a&\alpha a&0\\
b&c&\beta a
\end{array}\right]$&$
4,\,3$\\  & & & \\
7&$\left[\begin{array}{ccc}
e^{\alpha a}\cos(a)&e^{\alpha a}\sin(a)&b\\
-e^{\alpha a}\sin(a)&e^{\alpha a}\cos(a)&c\\
0&0&e^{\beta a}
\end{array}\right]$&$\left[\begin{array}{ccc}
\alpha a &a&b\\
-a&\alpha a&c\\
0&0&\beta a
\end{array}\right]$&$
4,\,3$\\  & & & \\
8&$\left[\begin{array}{ccc}
e^{\alpha a}\cos(a)&e^{\alpha a}\sin(a)&0\\
-e^{\alpha a}\sin(a)&e^{\alpha a}\cos(a)&0\\
0&0&e^{\beta a}
\end{array}\right]$&$\left[\begin{array}{ccc}
\alpha a &a&0\\
-a&\alpha a&0\\
0&0&\beta a
\end{array}\right]$&$
2,\,1$\\  & & & \\
9&$\left[\begin{array}{ccc}
a&b&c\\
d&e&f\\
g&h&i
\end{array}\right]$&$\left[\begin{array}{ccc}
a&b&c\\
d&e&f\\
g&h&i
\end{array}\right]$&$
8$\\  & & & \\
10&$\left[\begin{array}{ccc}
a&b&c\\
d&e&f\\
0&0&g
\end{array}\right]$&$\left[\begin{array}{ccc}
a&b&c\\
d&e&f\\
0&0&g
\end{array}\right]$&$
6$\\
& & & \\
\hline
\end{tabular}

\begin{tabular}{|cllc|}
\hline
Number&Lie subgroups with det$\,G=1$&Lie Subalgebras with tr$\,H=0$&Dimension\\
\hline
& & & \\
11&$\left[\begin{array}{ccc}
a&b&0\\
c&d&0\\
e&f&g
\end{array}\right]$&$\left[\begin{array}{ccc}
a&b&0\\
c&d&0\\
e&f&g
\end{array}\right]$&$
6$
\\  & & & \\
12&$\left[\begin{array}{ccc}
a&b&c\\
d&e&f\\
0&0&1
\end{array}\right]$&$\left[\begin{array}{ccc}
a&b&c\\
d&e&f\\
0&0&0
\end{array}\right]$&$
5$\\
13 &$\left[\begin{array}{ccc}
a&b&0\\
c&d&0\\
e& f&1
\end{array}\right]$&$\left[\begin{array}{ccc}
a&b&0\\
c&d&0\\
e&f&0
\end{array}\right]$&$
5$\\  & & & \\
14&$\left[\begin{array}{ccc}
a&b&0\\
c&d&0\\
0&0&e
\end{array}\right]$&$\left[\begin{array}{ccc}
a&b&0\\
c&d&0\\
0&0&e
\end{array}\right]$&$
4$\\  & & & \\
15&$\left[\begin{array}{ccc}
a&b&0\\
c&d&0\\
0&0&1
\end{array}\right]$&$\left[\begin{array}{ccc}
a&b&0\\
c&d&0\\
0&0&0
\end{array}\right]$&$
3$\\  & & & \\
16\footnotemark[1] & exp$\left[\begin{array}{ccc}
0&-a&-b\\
a&0&-c\\
b&c&0
\end{array}\right]$&$
\left[\begin{array}{ccc}
0&-a&-b\\
a&0&-c\\
b&c&0
\end{array}\right]$&$
3$\\ & & & \\
17&$\left[\begin{array}{ccc}
e^a&0&0\\
b&e^{-2a}&c\\
a e^a&0&e^a
\end{array}\right]$&$\left[\begin{array}{ccc}
a&0&0\\
b&-2a&c\\
a&0&a
\end{array}\right]$&$
3$\\  & & & \\
18&$\left[\begin{array}{ccc}
e^a&b&0\\
0&e^{-2a}&0\\
a e^a&c&e^a
\end{array}\right]$&$\left[\begin{array}{ccc}
a&b&0\\
0&-2a&0\\
a&c&a
\end{array}\right]$&$
3$\\  & & & \\
19\footnotemark[2] & exp$\left[\begin{array}{ccc}
0&a&b\\
b&c&0\\
a&0&-c
\end{array}\right]$&$
\left[\begin{array}{ccc}
0&a&b\\
b&c&0\\
a&0&-c
\end{array}\right]$&$
3$\\  & & & \\
20&$\left[\begin{array}{ccc}
1&0&a\\
a e^b&e^b&c\\
0&0&e^{-b}
\end{array}\right]$&$\left[\begin{array}{ccc}
0&0&a\\
a&b&c\\
0&0&-b
\end{array}\right]$&$
3$\\  & & & \\
\hline
\end{tabular}

\footnotetext[1]{The explicit form of this subgroup is very
complicated but is just the special orthogonal group.}

\footnotetext[2]{This subgroup is conjugate to the Lorentz group.}

\begin{tabular}{|cllc|}
\hline
Number&Lie subgroups with det$\,G=1$&Lie Subalgebras with tr$\,H=0$&Dimension\\
\hline
 & & & \\
21&$\left[\begin{array}{ccc}
e^a&0&0\\
b&e^{-2a}&0\\
a e^a&0&e^a
\end{array}\right]$&$\left[\begin{array}{ccc}
a&0&0\\
b&-2a&0\\
a&0&a
\end{array}\right]$&$
2$\\
& & & \\
22&$\left[\begin{array}{ccc}
e^a&0&0\\
0&e^{-2a}&0\\
a e^a&b&e^a
\end{array}\right]$&$\left[\begin{array}{ccc}
a&0&0\\
0&-2a&0\\
a&b&a
\end{array}\right]$&$
2$\\  & & & \\
23&$\left[\begin{array}{ccc}
1&0&a\\
a e^b&e^b&\frac{1}{2}a^2e^b\\
0&0&e^{-b}
\end{array}\right]$&$\left[\begin{array}{ccc}
0&0&a\\
a&b&0\\
0&0&-b
\end{array}\right]$&$
2$\\  & & & \\
24&$\left[\begin{array}{ccc}
1&a&0\\
0&1&0\\
a&b&1
\end{array}\right]$&$\left[\begin{array}{ccc}
0&a&0\\
0&0&0\\
a&b&0
\end{array}\right]$&$
2$\\  & & & \\
25&$\left[\begin{array}{ccc}
e^a&0&0\\
0&e^{-2a}&0\\
a e^a&0&e^a
\end{array}\right]$&$\left[\begin{array}{ccc}
a&0&0\\
0&-2a&0\\
a&0&a
\end{array}\right]$&$
1$\\  & & & \\
26&$\left[\begin{array}{ccc}
1&a&0\\
0&1&0\\
a&\frac{1}{2}a^2&1
\end{array}\right]$&$\left[\begin{array}{ccc}
0&a&0\\
0&0&0\\
a&0&0
\end{array}\right]$&$
1$\\
& & & \\
\hline
\end{tabular}
\end{center}

{\parindent 0cm

{\sc David Mar{\'\i}n}

{\it  Departament de Matem\`atiques,
Universitat Aut\`onoma de Barcelona,
08193 Bellaterra, Bar\-ce\-lo\-na, SPAIN\\
e-mail: davidmp@mat.uab.es}

\bigskip

{\sc Manuel de Le\'on}

{\it Instituto de Matem\'aticas y F{\'\i}sica Fundamental, CSIC,
Serrano 123, 28006 Madrid, Spain.\\
e-mail: mdeleon@imaff.cfmac.csic.es} }


\normalsize


\newpage

\bigskip


\begin{thebibliography}{99}



\bibitem{Bernard} D. Bernard: Sur la g\'eom\'etrie
diff\'erentielle des $G$-structures, {\sl Ann. Inst. Fourier} {\bf
10} (1960), 151-270.

\bibitem{Ber} D. Bertrand: {\sl Groupes alg\'{e}briques lin\'{e}aires
et th\'{e}orie de Galois diff\'{e}rentielle}, Cours de 3\`{e}me cycle,
1985-86, Universit\'{e} Paris VI.

\bibitem{Blair} D.E. Blair: {\sl Contact manifolds in Riemannian
geometry}, Springer-Verlag, Berlin, 1976.

\bibitem{Bloom} F. Bloom: {\sl Modern Differential Geometric
Techniques in the Theory of Continuous Distributions of
Dislocations}, Lecture Notes in Math. {\bf 733}, Springer, Berlin,
1979.

\bibitem{Borel} A. Borel: {\sl Linear Algebraic Groups}, Graduate
Texts in Mathematics {\bf 126}, Springer, New-York, 1991.


\bibitem{Chern} S.S. Chern: The geometry of $G$-structures,
{\sl Bull. Amer. Math. Soc.} {\bf 72} (1966), 167-219.

\bibitem{Che} C. Chevalley: {\sl La th\'{e}orie des groupes alg\'{e}briques},
Proceedings of the International Congress of Mathematicians, 14-21 August
1958 (Edimburg)

\bibitem{CDL} L.A. Cordero, C.T.J. Dodson,  M. de Le\'on:
{\sl Differential Geometry of Frame Bundles}, Mathematics and
Its Applications, Kluwer, Dordrecht, 1989.

\bibitem{EES} M. Elzanowski, M. Epstein, J. Sniatycki:
$G$-structures and material homogeneity, {\sl Journal of
Elasticity} {\bf 23} (1990), 167-180.

\bibitem{EdL2000} M. Epstein, M. de Le\'on: Homogeneity without
uniformity: towards a mathematical theory of functionally graded
materials, {\sl Int. J. of Solids Structures} {\bf 37} (2000),
7577-7591.


\bibitem{Fu} A. Fujimoto: {\sl Theory of $G$-structures},
Publications of the Study Group of Geometry, Vol. I. Tokyo,
1972.

\bibitem{Greub} W. Greub, S.Halperin, R. Vanstone:
{\sl Connections, Curvature, Cohomology}, Vol I. Acad. Press, New
York-London, 1972.


\bibitem{Kob2} E.T. Kobayashi: A remark on the existence of $G$-structure,
{\sl Proc. Amer. Math. Soc.} {\bf 16} (1959), 300-309.

\bibitem{Kob3} E.T. Kobayashi: A remark on the Nijenhuis tensor,
{\sl Pacific Journal of Math.} {\bf 12} (1962), 963-977.

\bibitem{KN} S. Kobayashi, K. Nomizu: {\sl Foundations of
Differential Geometry}, vol. I, Interscience Publishers, New
York, 1963.

\bibitem{Kob} S. Kobayashi: {\sl Transformations Groups in
Differential Geometry}, Springer, Berlin-New York, 1972.

\bibitem{LL} J. Lehmann-Lejeune: Sur l'integrabilite de certaines
G-structures, {\sl C.R. Acad. Sc. Paris} {\bf 258} (1964),
5326-5329.

\bibitem{LE1} M. de Le\'on, M. Epstein: On the integrability of second
order $G$-structures with applications to continuous theories of
dislocations. {\sl Reports on Mathematical Physics} {\bf 33} (3)
(1993), 419-436.

\bibitem{LE2} M. de Le\'on, M. Epstein: Corps mat\'eriels de degr\'e
sup\'erieur. {\sl Comptes Rendus Acad. Sc. Paris} {\bf 319},
S\'erie I, (1994), 615-620.

\bibitem{LE230} M. de Le\'on, M. Epstein:
Homogeneity conditions for generalized Cosserat media.
{\sl Journal of Elasticity {\bf 43} (1996), 189-201.}

\bibitem{LE231} M. de Le\'on, M. Epstein:
On uniformity of shells. {\sl Int. J. of Solids Structures}  {\bf
35} 17 (1998), 2173-2182.

\bibitem{LE3} M. de Le\'on, M. Epstein:
Geometrical Theory of Uniform Cosserat Media.
{\sl Journal of Geometry and Physics} {\bf 26} (1998), 127-170.
(Con M. Epstein).

\bibitem{Lie} S. Lie, F. Engel: {\sl Theorie der transformationsgruppen},
Vol 3. Leipzig. Teubner 1893.

\bibitem{Mack} K. Mackenzie: {\sl Lie groupoids and Lie algebroids in
Differential Geometry}, London Mathematical Society Lecture Note
Series {\bf 124}, Cambridge Univ. Press, Cambridge, 1987.

\bibitem{marsden1} J.E. Marsden, T.J.R. Hughes: {\sl Mathematical
Foundations of Elasticity}, Prentice Hall, New Jersey, 1983.

\bibitem{Maugin} G.A. Maugin: {\sl Material Inhomogeneities in
Elasticity}, London, Chapman \& Hall, 1993.

\bibitem{Noll} W. Noll: Materially Uniform Simple Bodies
with Inhomogeneities, {\sl Arch. Rational Mech. Anal.} {\bf 27},
(1967), 1-32.

\bibitem{Nono} T. N\^{o}no: On the symetry groups of simple materials:
An application of the theory of Lie groups, {\sl Journal Math.
Anal. Appl.} {\bf 24} (1968), 110-135.

\bibitem{Oni} A.L. Onishchik, E.B. Vinberg: {\sl Lie Groups
and Algebraic Groups}, Springer-Verlag, Berlin-Heidelberg, 1990.

\bibitem{Ramis} J.P. Ramis, J. Martinet, {\sl Th\'{e}orie de Galois
Diff\'{e}rentielle et Resommation}, Part 2 of {\sl Computer Algebra
and Differential Equations}, Academic Press, 1990.

\bibitem{TrToup} C. Truesdell, R.A. Toupin: {\sl Principles of
Classical Mechanics and Field Theory}, Handbuch der Physik, Vol.
III/1, Berlin-New York, Springer, 1960.

\bibitem{TrNoll} C. Truesdell, W. Noll: {\sl The Non-Linear
Field Theories of Mechanics}, Handbuch der Physik, Vol. III/3,
Berlin-New York, Springer, 1965.

\bibitem{WangTr} C.C. Wang, C. Truesdell: {\sl Introduction to
rational elasticity}, Noordhoff International Publishing,
Leyden, 1973.

\bibitem{Wang} C.C. Wang: On the Geometric Structures of
Simple Bodies, a Mathematical Foundation for the Theory of
Continuous Distributions of Dislocations, {\sl Arch. Rational
Mech. Anal.} {\bf 27}, (1967), 33-94.

\bibitem{Yano} K. Yano, M. Kon: {\sl Structures on manifolds}, World
Scientific, Singapore, 1984.


\end{thebibliography}
\end{document}